\documentclass[12pt]{smfart}
\usepackage{graphics}

\newtheorem{thm}{Th\'eor\`eme}[section]
\newtheorem*{theo}{Th\'eor\`eme}
\newtheorem{cor}[thm]{Corollaire}
\newtheorem{prop}[thm]{Proposition}
\newtheorem{lem}[thm]{Lemme}
\newtheorem{definition}[thm]{D\'efinition}
\newtheorem{remark}[thm]{Remarque}
\newtheorem{example}[thm]{Exemple}

\newtheorem{conj}{Conjecture}

\newtheorem*{lemme}{Lemme}

\author[Frank Loray]{Frank LORAY}
\address{(Charg\'e de Recherches au CNRS)\\
IRMAR, UFR de Math\'ematiques,\\ 
Universit\'e de Rennes 1, Campus de Beaulieu,\\ 
35042 Rennes Cedex (France)}
\email{frank.loray@univ-rennes1.fr}
\urladdr{http://name.math.univ-rennes1.fr/frank.loray}

\title{Sur les Th\'eor\`emes I et II de Painlev\'e}
\alttitle{Painlev\'e's Theorems I and II}

\begin{document}

\frontmatter
\begin{abstract}Un si\`ecle avant Alberto, Paul Painlev\'e \'etait lui aussi 
professeur \`a l'Uni\-ver\-sit\'e de Lille et posait les fondations 
des notions de feuilletage et d'holonomie \`a travers 
deux th\'eor\`emes que l'on trouve aujourd'hui au d\'ebut
des ``{\it Le\c cons de Stockholm}''. Je remercie Alberto
de m'avoir encourag\'e \`a les lire et je propose
ici d'en donner un \'enonc\'e pr\'ecis 
ainsi qu'une preuve rigoureuse illustr\'ee par de nombreux exemples.
Je termine par deux conjectures
quant au comportement global des applications d'holonomie, 
localement d\'efinies par le Th\'eor\`eme II.
\end{abstract}

\date{February 2004}
\dedicatory{en l'honneur d'Alberto Verjovsky}

\maketitle

\mainmatter

\section*{Introduction}\label{S:intro}

Paul Painlev\'e 
a implicitement introduit la notion de feuilletage holomorphe singulier
dans la premi\`ere partie de ses ``{\it Le\c cons de Stockholm}'' lorsqu'il \'etudie
les propri\'et\'es des \'equations diff\'erentielles alg\'ebriques du premier ordre
$F(y',y,x)=0$ (sp\'ecialement dans les ``{\it Le\c cons 2 et 3}'').

\eject

Afin de faciliter l'expos\'e, nous ne discu\-terons que des \'equations du premier degr\'e :
$$(E)\ \ \ \ \ {dy\over dx}={P(x,y)\over Q(x,y)}$$
avec $P,Q\in\mathbb C[x,y]$ polynomiaux dans chacune des variables $(x,y)\in\mathbb C^2$.
Si les conditions initiales $(x_0,y_0)\in\mathbb C^2$ satisfont $Q(x_0,y_0)\not=0$,
alors le th\'eor\`eme de Cauchy implique l'existence d'une unique solution $y=f(x)$
(pour l'\'equation $(E)$) qui est analytique sur un voisinage $D$ de $x_0$
et telle que $f(x_0)=y_0$. Nous noterons $f(x,x_0,y_0)$ cette solution.
Ce th\'eor\`eme est purement local.
Que peut-on dire des solutions globales~? Plus pr\'ecis\'ement,
que peut-on dire du prolongement analytique et des singularit\'es possibles
de $f(x,x_0,y_0)$ dans la variable $x$ ? ou encore $y_0$ ?

\section{Prolongement analytique et fonctions multiformes}

Rappelons les notions de prolongement analytique et de surface de Riemann d'une fonction multiforme.
Dans toute la suite, $\Omega$ d\'esigne un domaine (connexe) de $\mathbb C$,
$f$ un germe de fonction holomorphe en un point $x_0\in\Omega$,
$\gamma:[0,1]\to\Omega$ un chemin (continu ou diff\'erentiable) issu de $\gamma(0)=x_0$
et $x_1:=\gamma(1)$ son extr\'emit\'e.

\begin{definition}\rm
On dit qu'une suite de disques $D_0,D_1,\ldots,D_n\subset\Omega$
{\it recouvre} $\gamma$ si, pour un d\'ecoupage convenable $0=t_0<t_1<\cdots<t_n<t_{n+1}=1$
de l'intervalle, $D_k$ contient $\gamma([t_k,t_{k+1}])$ pour $k=0,\ldots,n$.
On dit que $f$ admet un {\it prolongement analytique le long de} $\gamma$
s'il existe une suite de disque $D_0,D_1,\ldots,D_n\subset\Omega$ recouvrant $\gamma$
et des fonctions holomorphes $f_k:D_k\to\mathbb C$ telles que
$f_0\equiv f$ au voisinage de $x_0$ et $f_k\equiv f_{k-1}$
sur $D_{k-1}\cap D_k$ pour tout $k=1,\ldots,n$.
Le germe de fonction d\'efini par $f_n$ en $x_1:=\gamma(1)$
sera not\'e $f_\gamma$ et appel\'e {\it d\'etermination de $f$ au dessus de $x_1$}.
\end{definition}

\begin{remark}\rm
Dans la d\'efinition pr\'ec\'edente, la d\'etermination $f_\gamma$ ne d\'epend
que de la chaine de disques $(D_k)_k$.
En particulier, pour $\varepsilon>0$ suffisamment petit,
$f$ admet un prolongement analytique le long
de toute {\it $\varepsilon$-perturbation} de $\gamma$,
c'est \`a dire de tout chemin $\gamma':[0,1]\to\Omega$ satisfaisant
$\vert\gamma'(t)-\gamma(t)\vert<\varepsilon$ et $\gamma'(0)=x_0$.
De plus, d\`es que $\gamma'(1)=x_1$, ce nouveau chemin nous conduit
\`a la m\^eme d\'etermina\-tion $f_{\gamma'}\equiv f_\gamma$.

\begin{figure}[htbp]
\begin{center}

\input{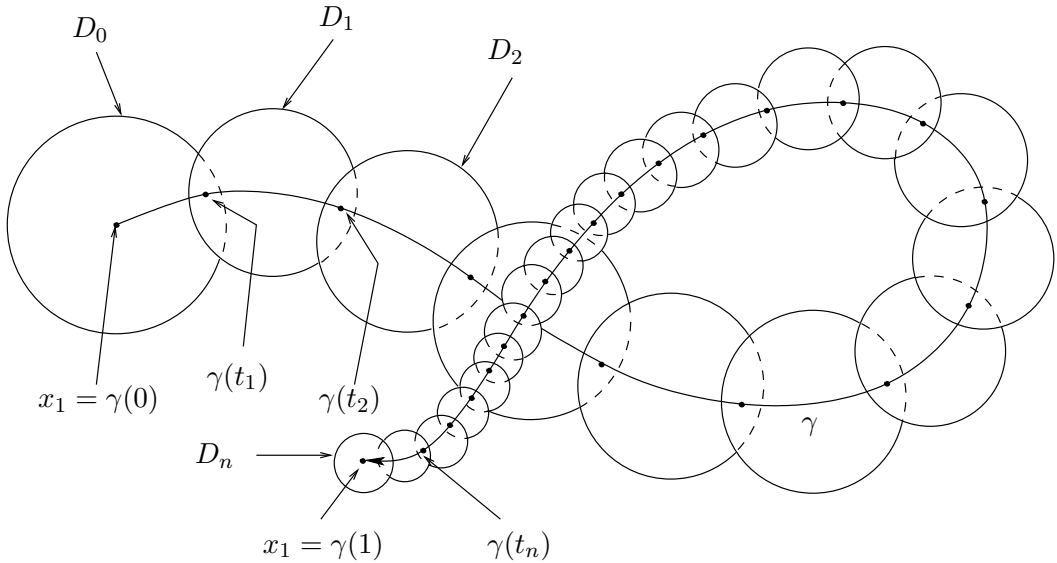}
 
\caption{Prolongement analytique}
\label{figure:1}
\end{center}
\end{figure}

\noindent Ainsi, quitte \`a remplacer $\gamma$ par une perturbation,
on pourra toujours supposer $\gamma$ diff\'erentiable dans la d\'efinition pr\'ec\'edente
et celles qui suivent.
\end{remark}

\begin{lemme}[Poincar\'e-Volterra]
Au dessus de chaque point $x_1\in\Omega$, l'ensemble :
$$\{f_\gamma\ ;\ \gamma:[0,1]\to\Omega,\ \gamma(0)=x_0\ \text{et}\ \gamma(1)=x_1\}$$
des d\'eterminations de $f$ est au plus d\'enombrable.
\end{lemme}

\begin{proof}
Toute d\'etermination $f_\gamma$ est d\'etermin\'ee
par la suite de disques $(D_k)_k$ recouvrant $\gamma$.
Quitte \`a remplacer les $D_k$ par des disques tr\`es proches
(ce qui ne modifie pas le germe de fonction obtenu en $x_1$)
on peut supposer leur rayon et les coordonn\'ees de leur centre rationnels.
L'ensemble des suites finies de disques \`a coordonn\'ees rationnelles 
\'etant d\'enombrable,
l'ensemble des d\'eterminations l'est aussi.
\end{proof}

\eject

\begin{prop}
Il existe une surface de Riemann (connexe) $\mathcal S$ avec un point marqu\'e $p_0\in\mathcal S$
et une application holomorphe $\phi:\mathcal S\to\Omega\times\mathbb C$, $\phi=(\pi,\widetilde f)$,
satisfaisant~:
\begin{enumerate}
\item $\pi(p_0)=x_0$ et $f\circ\pi\equiv\widetilde f$ au voisinage de $p_0$ ;
\item $\pi:\mathcal S\to\Omega$ est un diff\'eomorphisme local en tout point de $\mathcal S$ ;
\item tout autre triplet $(\mathcal S',p_0',\phi')$ satisfaisant (1) et (2)
se factorise par $(\mathcal S,p_0,\phi)$ via une application holomorphe $\varphi:\mathcal S'\to\mathcal S$ :
$\varphi(p_0')=p_0$ et $\phi'=\phi\circ\varphi$.
\end{enumerate}
Le triplet $(\mathcal S,p_0,\phi)$ est unique \`a isomorphisme pr\`es.
Le {\it graphe} $\mathcal G:=\phi(\mathcal S)$ est unique.
\end{prop}

\begin{figure}[htbp]
\begin{center}

\input{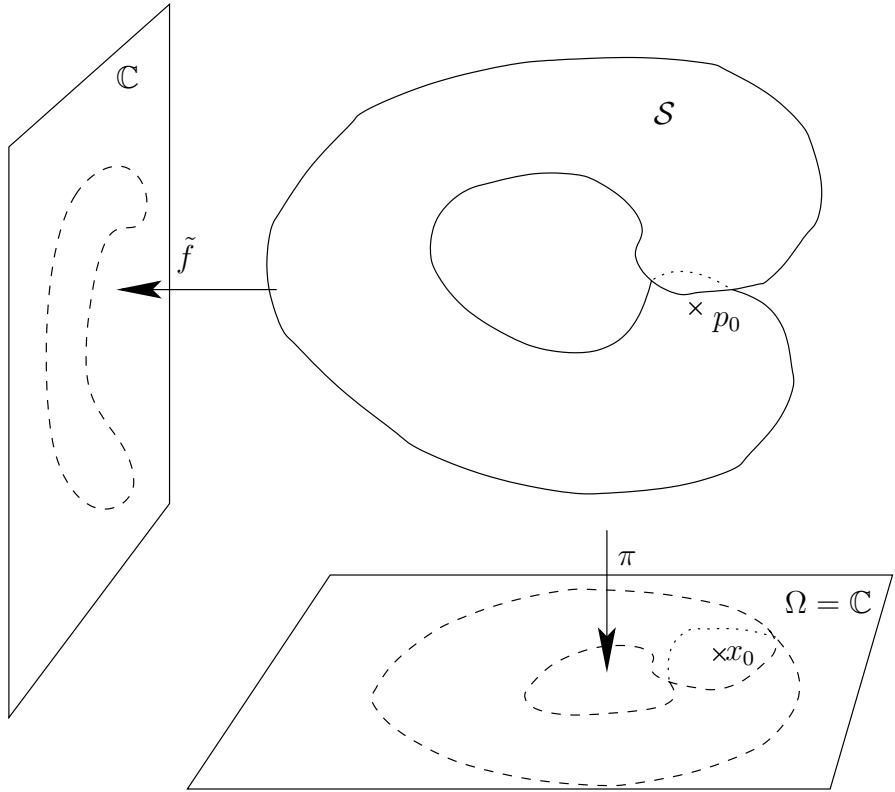}
 
\caption{Surface de Riemann (ou graphe) de $f$}
\label{figure:2}
\end{center}
\end{figure}

\begin{proof}[Id\'ee de preuve (voir \cite{Va})]
La surface est construite comme l'ensemble des couples
$(\gamma(1),f_\gamma)$ o\`u $\gamma$ d\'ecrit l'ensemble des chemins $\gamma:[0,1]\to\Omega$
issus de $\gamma(0):=x_0$ le long desquels $f$ admet un prolongement analytique.
Si $D$ est le disque de convergence de $f_\gamma$,
alors l'ensemble des couples $(\gamma'(1),f_{\gamma'})$ obtenus par concat\'enations
$\gamma'=\gamma''\cdot\gamma$ avec les chemins $\gamma'':[0,1]\to D$
issus de $\gamma''(0):=x_1$ forme un voisinage $V(p_1)$ du point $p_1:=(x_1,f_\gamma)$
de la surface. Par Poincar\'e-Volterra, une collection d\'enombrable de telles cartes
suffit pour recouvrir la surface. L'application $\phi$ est alors d\'efinie sur $V(p_1)$
par la fl\`eche $(x_1',f_{\gamma'})\mapsto(x_1',f_{\gamma'}(x_1'))$.
\end{proof}

\begin{remark}\rm
Le germe $f$ admet un prolongement analytique le long de $\gamma$
si et seulement si $\gamma$ se rel\`eve en $\widetilde\gamma:[0,1]\to\mathcal S$ satisfaisant
$\widetilde\gamma(0)=p_0$ et $\gamma\equiv\pi\circ\widetilde\gamma$.
Alors, $f_\gamma$ est le germe de fonction d\'efini en $\gamma(1)$ par
l'\'egalit\'e $f_\gamma\circ\pi=\widetilde f$ en $\widetilde\gamma(1)$.
\end{remark}

\begin{definition}\rm
On appellera {\it fonction multiforme d\'efinie sur $\Omega$ par le germe $f$}
le ``prolongement analytique maximal'' de $f$ sur $\Omega$ donn\'e
ou bien par sa surface de Riemann $\phi:\mathcal S\to\Omega\times\mathbb C$
(\`a isomorphisme pr\`es), ou bien par son graphe $\mathcal G$.
On notera $\underline f_{\Omega}$ cette ``fonction'' et, si n\'ecessaire,
$\phi_{\underline f_{\Omega}}:\mathcal S_{\underline f_{\Omega}}\to\Omega\times\mathbb C$
et $\mathcal G_{\underline f_{\Omega}}$.
\end{definition}

\begin{remark}\rm La restriction d'une fonction multiforme sur $\Omega$
\`a un domaine $\Omega'\subset\Omega$ n'a pas de sens :
la pr\'eimage $\phi^{-1}_{\underline f_{\Omega}}(\Omega'\times\mathbb C)$ peut avoir plusieurs
composantes connexes dans $\mathcal S_{\underline f_{\Omega}}$ et d\'efinir autant de fonctions multiformes
sur $\Omega'$. Par contre, si $x_0\in\Omega'$, alors la composante connexe contenant $p_0$
s'identifie \`a $\mathcal S_{\underline f_{\Omega'}}$.
\end{remark}

\section{Singularit\'es}

\begin{definition}\rm
On dit que $\gamma$ {\it conduit $f$ vers une singularit\'e}
lorsque $f$ admet un prolongement analytique le long de
$\gamma\vert_{[0,1-\varepsilon]}$ pour tout $\varepsilon>0$
mais pas le long de $\gamma$.
On dit qu'un autre chemin $\gamma':[0,1]\to\Omega$ joignant $x_0$ \`a $x_1$
conduit $f$ vers {\it la m\^eme singularit\'e} si pour tout voisinage $D$ de $x_1$,
les d\'eterminations $f_{\gamma\vert_{[0,\tau]}}$
et $f_{\gamma'\vert_{[0,\tau]}}$ d\'efinissent la m\^eme fonction multiforme sur $D$
pour $\tau$ suffisamment proche de $1$.

On appelle {\it singularit\'e de $f$ au dessus de $\Omega$} (ou de $\underline f_\Omega$)
{\it d\'efinie par $\gamma$} la classe d'\'equivalence not\'ee $f_\gamma$ pour la
relation pr\'ec\'edente.
On note $\widetilde\Sigma_{\underline f_\Omega}$ (resp. $\Sigma_{\underline f_\Omega}$)
l'ensemble des singularit\'es $f$ au dessus de $\Omega$
(resp. de leur projection sur $\Omega$).
\end{definition}

\begin{figure}[htbp]
\begin{center}

\input{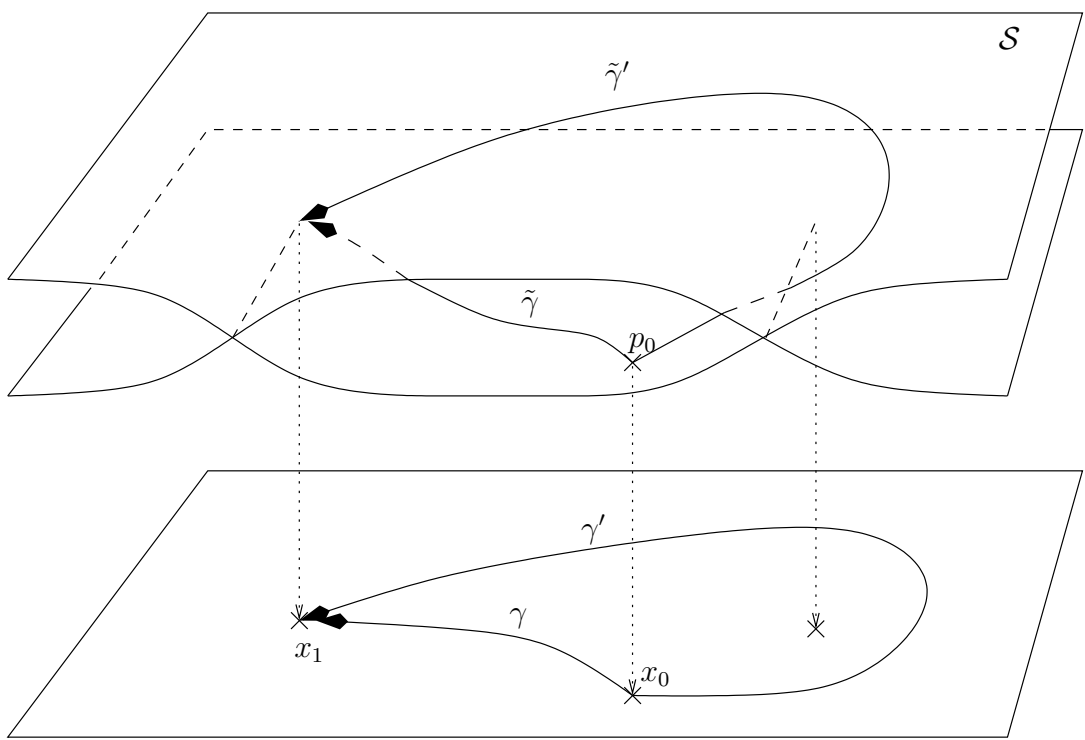}
 
\caption{Deux chemins conduisant vers une m\^eme singularit\'e}
\label{figure:3}
\end{center}
\end{figure}

\begin{remark}\rm
En fait, $\gamma$ conduit $f$ vers une singularit\'e si et seulement si
le chemin incomplet $\gamma\vert_{[0,1[}:[0,1[\to\Omega$
se rel\`eve, via la projection $\pi:\mathcal S\to\Omega$,
en un chemin propre $\widetilde\gamma:[0,1[\to\mathcal S$
(c'est \`a dire $\widetilde\gamma(t)$ tend
vers le bord de $\mathcal S$ lorsque $t\to 1$).
\end{remark}

\eject

\begin{definition}\rm
On dit que $f$ d\'efinit une {\it fonction multiforme r\'eguli\`ere} sur $\Omega$
si $f$ ne poss\`ede pas de singularit\'e au dessus de $\Omega$, i.e. admet un prolongement analytique le long de tout chemin $\gamma:[0,1]\to\Omega$
issu de $\gamma(0)=x_0$. C'est le cas si et seulement si l'application
$\pi:\mathcal S\to\Omega$ est un rev\^etement. En particulier,
le nombre de d\'eterminations de $f$ au dessus d'un point $x_1\in\Omega$
est independant de $x_1$ ; si $\Omega$ est simplement connexe,
alors l'extension est uniforme.
\end{definition}

\begin{example}\label{E:logarithme}\rm
La d\'etermination principale du logarithme $f(x)=\log(x)$ en $x_0=1$,
$f(1)=0$, d\'efinit une fonction uniforme sur $\mathbb C\setminus\mathbb R^-$,
multiforme r\'eguli\`ere sur $\mathbb C^*$
et poss\`ede, au dessus de $\mathbb C$, une unique singularit\'e $\gamma$
qui se projette sur $0$.
\end{example}

\begin{example}\label{E:puissanceComplexe}\rm
La fonction multiforme d\'efinie par $f(x)=x^\alpha:=\exp(\alpha\log(x))$,
$\alpha\in\mathbb C^*$, solution de l'\'equation diff\'erentielle $y'=\alpha{y\over x}$,
poss\`ede exactement une singularit\'e au dessus de $0$
except\'e lorsque $\alpha\in\mathbb N$.
Elle poss\`ede une infinit\'e de d\'eterminations au dessus de chaque point autre que $0$
sauf quand $\alpha$ est rationnel~:
si $\alpha={p\over q}\in\mathbb Q$, alors $f$ poss\`ede exactement $q$ d\'eterminations.
Enfin, $f$ admet une limite en $0$ si et seulement si $\alpha$ est r\'eel :
$f(x)$ tend vers $0$ si $\alpha>0$ et $\infty$ si $\alpha<0$
d\`es que $x\to 0$ le long de n'importe quel chemin $\gamma$ repr\'esentant la singularit\'e.
\end{example}

\begin{example}\label{E:RiccatiRetournee}\rm
La fonction
$f(x)=\alpha_1\log(x-\zeta_1)+\alpha_2\log(x-\zeta_2)+\alpha_3\log(x-\zeta_3)$,
solution de l'\'equation diff\'erentielle
${dy\over dx}={\alpha_1\over x-\zeta_1}+{\alpha_2\over x-\zeta_2}+{\alpha_3\over x-\zeta_3}$,
est r\'eguli\`ere sur $\mathbb C\setminus\{\zeta_1,\zeta_2,\zeta_3\}$.
Au dessus d'un point distinct des $\zeta_k$, deux d\'eterminations diff\`erent d'une constante
appartenant au groupe additif $G=2\pi i(\alpha_1\mathbb Z+\alpha_2\mathbb Z+\alpha_3\mathbb Z)$.
Au dessus de $\zeta_1$, deux singularit\'es diff\`erent d'une constante
appartenant au sous-groupe $H_1=2\pi i(\alpha_2\mathbb Z+\alpha_3\mathbb Z)$.
La fonction inverse $x=g(y)$, solution de l'\'equation diff\'erentielle
${dx\over dy}=1/
\left({\alpha_1\over x-\zeta_1}+{\alpha_2\over x-\zeta_2}+{\alpha_3\over x-\zeta_3}\right)$,
poss\`ede au moins comme singularit\'es les valeurs critiques de $f$.
En g\'en\'eral,
la d\'eriv\'ee $f'(x)$ s'annule au dessus de deux points $x_1\not=x_2$ distincts des $\zeta_k$
et pour toutes les d\'etermination. Pour des $\alpha_k$ g\'en\'eriques,
le groupe $G$ est dense dans $\mathbb C$ ; par suite, les singularit\'es de $g$
se projettent sur un sous-ensemble $\Sigma$ dense de $\mathbb C$.
\end{example}

\eject

\begin{figure}[htbp]
\begin{center}

\input{Alberto4.pstex_t}
 
\caption{Exemple \ref{E:RiccatiRetournee}}
\label{figure:4}
\end{center}
\end{figure}

\begin{example}\label{E:SingNonIsolee}\rm
La fonction multiforme d\'efinie sur $\mathbb C$ par $f(x)=\sqrt{x^\alpha-1}$,
$\alpha={\log(2)\over 2\pi i}$,
solution de l'\'equation diff\'erentielle $y'={\alpha(y^2+1)\over 2xy}$,
poss\`ede exactement une singularit\'e au dessus de chaque point $x_n=2^n$ et au dessus de $0$.
Pr\'ecis\'ement, sur $\mathbb C\setminus\mathbb R^-$ et pour la d\'etermination principale
de $x^\alpha=\exp(\alpha\log(x))$,
$f$ poss\`ede exactement une singularit\'e au dessus de $x_0=1$ et $2$ d\'eterminations autour :
notons $f^+$ et $f^-$ les fonctions uniformes induites sur
$\mathbb C\setminus]-\infty,0]\cup[1,+\infty[$.
Apr\`es avoir tourn\'e $n$ fois autour de $0$, $n\in\mathbb Z$,
la situation est essentiellement la m\^eme au dessus de $\mathbb C\setminus\mathbb R^-$ :
l'unique singularit\'e est maintenant au dessus $x_n=2^n$ et les $2$ nouvelles
d\'eterminations uniformes sur $\mathbb C\setminus]-\infty,0]\cup[x_n,+\infty[$
respectivement not\'ees $f_n^+$ et $f_n^-$ se permuttent autour de $x_n$.
Ainsi, au dessus de chacun des $x_n$, il y a exactement une singularit\'e
et une infinit\'e de d\'eterminations r\'eguli\`eres.
Au dessus de $0$, tous les chemins m\^enent \`a la m\^eme singularit\'e :
\`a l'int\'erieur de tout disque $D$ centr\'e en $0$, toutes les d\'eterminations
$f_n^+$ et $f_n^-$ (en un point g\'en\'erique) s'\'echangent
par prolongement analytique dans $D$. En effet, on passe de $f_n^+$ \`a $f_{n+1}^+$
(resp. de $f_n^-$ \`a $f_{n+1}^-$) en tournant autour de $0$ ; enfin,
pour $n>>0$, $x_n$ est dans $D$ et on passe de $f_n^+$ \`a $f_n^-$ en tournant autour.
Cette singularit\'e est plus compliqu\'ee que les pr\'ec\'edentes
car d'autres singularit\'es l'accumulent.
\end{example}

\eject

\begin{figure}[htbp]
\begin{center}

\input{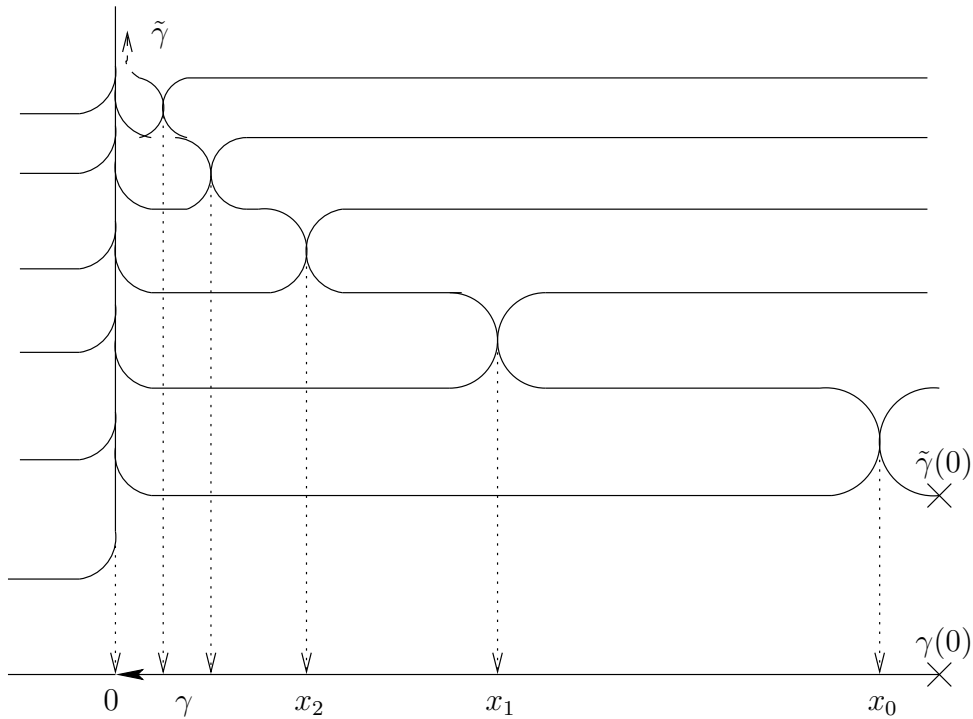}
 
\caption{La singularit\'e non isol\'ee $\tilde\gamma$ de l'exemple \ref{E:SingNonIsolee}}
\label{figure:5}
\end{center}
\end{figure}

\begin{example}\label{E:Picard}\rm
Le rev\^etement universel
$\{\vert x\vert<1\}\to\mathbb C\setminus\{0,1\}$
est une fonction uniforme avec fronti\`ere naturelle~:
il n'existe pas de prolongement analytique au del\`a du cercle
$\Sigma=\{\vert x\vert=1\}$.
\end{example}

Les diff\'erents types de singularit\'es apparaissant dans les exemples pr\'ec\'e\-dents
se distinguent modulo quelques d\'efinitions suppl\'ementaires.

\begin{definition}\rm
Une singularit\'e $f_\gamma$ est dite {\it isol\'ee}
s'il existe un disque $D$ centr\'e en $x_1:=\gamma(1)$ tel que
la fonction multiforme $\underline{f_\gamma}_D$ qu'elle d\'efinit sur $D$
soit r\'eguli\`ere sur $D^*:=D\setminus\{x_1\}$.
\end{definition}

\begin{remark}\rm
La singularit\'e d\'efinie par $\gamma$ est isol\'ee si et seulement si
il existe $\varepsilon>0$ tel que $f$ admette un prolongement analytique
le long de tout chemin $\gamma':[0,1]\to\Omega$ issu de $x_0$,
$\varepsilon$-proche de $\gamma$ \'evitant $x_1$~:
$$\vert\gamma'(t)-\gamma(t)\vert<\varepsilon\ \ \ \text{et}\ \ \ x_1\not\in\gamma'(]0,1]).$$
\end{remark}

\begin{definition}\rm
On dit que $f$ admet $y_1\in\overline{\mathbb C}$ comme limite en une singularit\'e $f_\gamma$,
$\overline{\mathbb C}=\mathbb C\cup\infty$,
si pour tout chemin $\gamma'$ repr\'esen\-tant cette singularit\'e
on a $\lim_{t\to1}f(\gamma'(1))=y_1$.

Une singularit\'e $f_\gamma$ est dite {\it alg\'ebro\"\i de}
si elle est isol\'ee, poss\`ede un nombre fini de d\'eterminations
autour de $x_1:=\gamma(1)$ et poss\`ede une limite (finie ou infinie) en $x_1$.
\end{definition}

\begin{remark}\rm
Soit $f_\gamma$ une singularit\'e isol\'ee. Sont \'equivalents :
\begin{enumerate}
\item $f_\gamma$ poss\`ede un nombre fini de d\'eterminations autour de $x_1=\gamma(1)$,
\item le relev\'e $\widetilde\gamma:[0,1[\to\mathcal S$ de $\gamma$,
tend vers une composante parabolique du bord de $\mathcal S$
(i.e. vers un point $p_1$ dans une surface de Riemann $\widehat{\mathcal S}=\mathcal S\cup\{p_1\}$),
\item le graphe de la fonction multiforme $f_\gamma$ sur un disque $D$
suffisamment petit centr\'e en $x_1$
est un disque \'epoint\'e proprement plong\'e dans $D^*\times\mathbb C$,
\item le germe de fonction multiforme $f_\gamma$ est de la forme
$g\left((x-x_1)^{1/q}\right)$
o\`u $g$ est une fonction holomorphe sur un voisinage \'epoint\'e de $0\in\mathbb C$ et $q\in\mathbb N^*$.
\end{enumerate}
On dit alors que {\it la singularit\'e est \`a monodromie finie} 
et sont \'equivalents~:
\begin{enumerate}
\item $f_\gamma$ est alg\'ebro\"\i de,
\item $\widetilde f$ s'\'etend m\'eromorphiquement en $p_1$,
\item $\widetilde{D^*}$ est la restriction \`a $D^*\times\mathbb C$
d'un disque $\widetilde D\subset D\times\overline{\mathbb C}$
proprement plong\'e,
\item $g$ s'\'etend m\'eromorphiquement en $0$.
\end{enumerate}
\end{remark}

Les singularit\'es alg\'ebro\"\i des g\'en\'eralisent
l'exemple \ref{E:puissanceComplexe} avec $\alpha$ rationel~;
ce sont aussi les singularit\'es des fonctions alg\'ebriques.
Dans l'exemple \ref{E:RiccatiRetournee}, toutes les singularit\'es 
de $g$ au dessus de $\mathbb C$ sont alg\'ebro\"\i des
(et donc isol\'ees) bien que leur projection soit dense.
Dans l'exemple \ref{E:SingNonIsolee}, toutes les singularit\'es 
sont alg\'ebro\"\i des
except\'ee celle au dessus de $0$ qui n'est m\^eme pas isol\'ee.

\eject

\section{Fonctions alg\'ebro\"\i des}

\begin{definition}\rm
On dit que $f$ d\'efinit une {\it fonction alg\'ebro\"\i de} sur $\Omega$
si toute singularit\'e de $f$ au dessus de $\Omega$ est alg\'ebro\"\i de.
\end{definition}

\begin{example}\rm
Si $\pi:\mathcal S\to\Omega$ est une application holomorphe
propre (non n\'ecessai\-rement finie) d'une surface de Riemann $\mathcal S$ sur $\Omega$
et si $g:\mathcal S\to\mathbb C$ est une fonction m\'eromorphe,
alors la fonction multiforme $f$ d\'efinie sur $\Omega$ par
$(\pi,g):\mathcal S\to\Omega\times\mathbb C$ est alg\'ebro\"\i de sur $\Omega$.
Les fonctions alg\'ebro\"\i des ainsi construites satisfont en outre \`a la propri\'et\'e suivante.
Elles ne poss\`edent qu'un nombre fini de d\'eterminations d\`es que l'on restreint
le prolongement analytique \`a un ouvert relativement compact dans $\Omega$.
La fonction $f$ de l'exemple \ref{E:RiccatiRetournee} est alg\'ebro\"\i de 
sur $\mathbb C$ mais ne satisfait pas
\`a cette derni\`ere propri\'et\'e.\end{example}

\begin{prop}\label{P:ProlongementAlgebroide}
Si $f$ d\'efinit une fonction alg\'ebro\"\i de sur $\Omega$, alors :
\begin{enumerate}
\item l'ensemble $\widetilde\Sigma$ des singularit\'es de $f$
au dessus de $\Omega$ est d\'enombrable,
\item pour tout chemin $\gamma:[0,1]\to\Omega$ issu de $x_0$,
il existe un sous-ensemble fini $\Sigma(\gamma)\subset\Sigma\cap\gamma(]0,1])\subset\Omega$
et un $\varepsilon>0$ tels que $f$ se prolonge analytiquement le long de tout chemin
$\varepsilon$-proche de $\gamma$ \'evitant $\Sigma(\gamma)$ :

\noindent $\gamma':[0,1]\to\Omega$,\hfill $\gamma(0)=x_0$,\hfill
$\vert\gamma'(t)-\gamma(t)\vert<\varepsilon$\hfill
et\hfill $\gamma'(]0,1[)\cap\Sigma(\gamma)=\emptyset$

\noindent conduisant \'eventuellement $f$ vers une singularit\'e
lorsque $\gamma'(1)=\gamma(1)$.
\end{enumerate}
En particulier, l'ensemble des valeurs prises par les diff\'erentes 
d\'etermi\-nations
$\{f_{\gamma'}(x_1)\ ;\ \gamma'\ \text{comme au dessus avec}\
\gamma'(1)=x_1\}\subset\overline{\mathbb C}$
est fini.
\end{prop}

\begin{remark}\rm
Dans la proposition pr\'ec\'edente, l'ensemble $\Sigma(\gamma)$
joue le r\^ole des singularit\'es interm\'ediaires possibles qu'il suffit de contourner
(en perturbant $\gamma$)
de fa\c con \`a poursuivre le prolongement analytique de $f$ jusqu'\`a $x_1$.
\end{remark}

Cette proposition est d\'emontr\'ee au d\'ebut de \cite{We} ; la preuve 
que nous donnons ici a l'avantage de montrer qu'elle reste vraie,
\`a la finitude pr\`es de $\Sigma(\gamma)$, sous l'hypoth\`ese
beaucoup plus faible que toute singularit\'e de $f$ sur $\Omega$ est isol\'ee.

\eject

\begin{lem}
Si toute singularit\'e de $f$ au dessus de $\Omega$ est isol\'ee,
alors l'ensemble $\widetilde\Sigma$ des singularit\'es de $f$
au dessus de $\Omega$ est au plus d\'enombrable.
\end{lem}

\begin{proof}
C'est une adaptation du lemme de Poincar\'e-Volterra.
\`A toute singularit\'e $f_\gamma$, on peut associer une suite de disques
$D_0,D_1,\ldots,D_n\subset\Omega$ recouvrant $\gamma$
et des fonctions holomorphes $f_k:D_k\to\mathbb C$ pour $k=1,\ldots,n-1$ tels que
$f_0\equiv f$ au voisinage de $x_0$, $f_k\equiv f_{k-1}$ sur $D_{k-1}\cap D_k$
et la fonction multiforme d\'efinie sur $D_n$ par $f_{n-1}$
poss\`ede exactement une singularit\'e au dessus de $D_n$, \`a savoir $f_\gamma$.
Quitte \`a diminuer la taille des disques,
nous pouvons de plus supposer leur rayon
et les coordonn\'ees de leur centre dans $\mathbb Q$.
\'Etant donn\'ee une telle suite de disques $D_k$, les fonctions $f_k$
sont bien d\'efinies par unicit\'e du prolongement analytique ;
il en est de m\^eme de la singularit\'e $f_\gamma$ au dessus de $x_1$ :
on peut oublier le chemin.
Puisque l'ensemble de telles suites de disques est d\'enombrable, l'ensemble des
singularit\'es l'est aussi.
\end{proof}

\begin{remark}\label{R:Contournement}\rm
Supposons que l'on rencontre une singularit\'e de $f$ en $t_1<1$
le long d'un chemin $\gamma:[0,1]\to \Omega$ issu de $x_0=\gamma(0)$
(i.e. le chemin $\gamma\vert_{[0,t_1]}$ conduit $f$
vers une singularit\'e) ; si cette singularit\'e est isol\'ee,
alors on peut la contourner dans le sens suivant.
Pour $\varepsilon>0$ suffisamment petit,
on peut modifier $\gamma$ le long de la composante connexe $]t_1^-,t_1^+[$ de
$\gamma^{-1}(\{\vert x-x_1\vert<\varepsilon\})$
contenant $t_1$ de sorte que le nouveau chemin $\gamma_1$
obtenu ainsi \'evite $x_1:=\gamma(t_1)$ :
$\gamma(]t_1^-,t_1^+[)\subset\{\vert x-x_1\vert<\varepsilon\}\setminus\{x_1\}$.
Le long de $\gamma_1$, le prolongement analytique de $f$ peut \^etre poursuivi
un peu plus loin, i.e. au moins jusqu'\`a $t_1^+$. Attention,
la d\'etermination de $f$ obtenue en $\gamma(t_1^+)$ d\'epend
de la perturbation $\gamma_1$ choisie ;
cependant, parce que la singularit\'e est isol\'ee,
elle ne d\'epend que de la classe d'homotopie de $\gamma_1$ parmi
ceux donn\'es par la construction pr\'ec\'edente.
Pour une singularit\'e non isol\'ee, cette construction ne serait plus
v\'eritablement locale et perdrait de son sens.
\end{remark}

\eject

\begin{figure}[htbp]
\begin{center}

\input{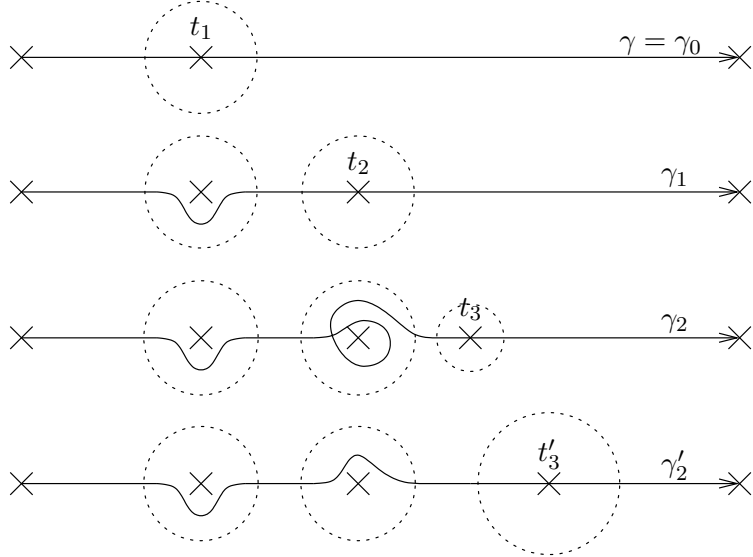}
 
\caption{Contournement successif des singularit\'es}
\label{figure:6}
\end{center}
\end{figure}

\begin{lem}
Si toute singularit\'e de $f$ au dessus de $\Omega$ est isol\'ee,
alors pour tout chemin $\gamma:[0,1]\to \Omega$
issu de $x_0=\gamma(0)$,
quitte \`a r\'ep\'eter la construction de la remarque pr\'ec\'edente
un nombre fini de fois,
on obtient un chemin $\gamma_n$ le long duquel $f$ admet un prolongement
analytique
(ou conduit vers une singularit\'e).
\end{lem}

En d'autres termes, on n'aura toujours qu'un nombre fini de singularit\'es
\`a contourner avant d'atteindre $\gamma(1)$,
et ceci ind\'ependamment de la fa\c con dont on
contourne ces ``singularit\'es interm\'ediaires''.

\begin{proof}
Supposons que pour certains choix successifs dans la fa\c con
de contourner les singularit\'es interm\'ediaires,
nous soyons amen\'es \`a effectuer une infinit\'e de contournements successifs
 $\gamma_0=\gamma,\gamma_1,\gamma_2,\ldots$
mettant en \'evidence une infinit\'e de singularit\'es successives $0<t_1<t_2<\cdots<1$,
avec limite $0<t_\infty\le1$. Alors on peut choisir les modifications successives
$\gamma_n\to\gamma_{n+1}$ de plus en plus petites de sorte que la suite
$\gamma_n$ d'applications continues tende uniform\'ement vers un chemin
$\gamma_\infty$.
Alors, $f$ poss\`ede une singularit\'e en $\gamma(t_\infty)$
pour le prolongement analytique
le long de $\gamma_\infty:[0,t_\infty]\to\Omega$ qui n'est, par construction,
 pas isol\'ee.
\end{proof}

\begin{remark}\rm
Les singularit\'es successives $0<t_1<t_2<\ldots<t_n<1$
rencontr\'ees durant l'algorithme du Lemme pr\'ec\'edent
ainsi que leur nombre $n$ d\'ependent des choix successifs que l'on fait
en les contournant.
Par exemple, si la premi\`ere singularit\'e $t_1$
(qui, elle, ne d\'epend d'aucun choix) n'est pas \`a monodromie finie,
alors il existe une infinit\'e de fa\c cons de la contourner
et la position de la singularit\'e suivante $t_2$ peut d\'ependre
enti\`erement de ce choix.
\end{remark}

\begin{lem}
Si toute singularit\'e de $f$ au dessus de $\Omega$ est isol\'ee
\`a monodromie finie, alors pour tout
chemin $\gamma:[0,1]\to\Omega$ issu de $\gamma(0)=x_0$,
il existe un sous-ensemble fini
$\Sigma(\gamma)\subset\Sigma\cap\gamma(]0,1[)\subset\Omega$
tel que l'on ait la propri\'et\'e (2) de la proposition
\ref{P:ProlongementAlgebroide}.
\end{lem}

\begin{proof}
Supposons qu'il existe une infinit\'e de singularit\'es interm\'ediaires
pour le prolongement alg\'ebroide le long de $\gamma$. Puisque la premi\`ere singularit\'e
$\gamma(t_1)$ est \`a monodromie finie, il n'existe qu'un nombre fini de fa\c cons
essentiellement distinctes de la contourner.
Pour au moins un de ces choix $\gamma_1$, il doit rester une infinit\'e de singularit\'es
interm\'ediaires possibles pour le prolongement analytique le long de $\gamma_1$.
En r\'ep\'etant cette op\'eration ind\'efiniment, on est amen\'e \`a contourner
une infinit\'e de singularit\'es le long de $\gamma$,
ce qui contredit le Lemme pr\'ec\'edent.
\end{proof}

\section{Les th\'eor\`emes de Painlev\'e :
des \'equations diff\'erentielles aux feuilletages}

Nous sommes maintenant en mesure d'\'enoncer le premier th\'eor\`eme 
de Painlev\'e~:

\begin{theo}[I]
Soient $P,Q\in\mathbb C[x,y]$
et consid\'erons l'\'equation diff\'erentielle :
$$(E)\ \ \ \ \ {dy\over dx}={P(x,y)\over Q(x,y)}.$$
Il existe un ensemble fini de points $\Sigma_E\subset\mathbb C$ 
sur lequel se projettent toutes les singularit\'es
non alg\'ebro\"\i des au dessus de $\mathbb C$ de toute solution
locale $y=f(x,x_0,y_0)$ de $(E)$.
\end{theo}

Ainsi, toute solution est une fonction alg\'ebro\"\i de
sur $\Omega=\mathbb C\setminus\Sigma_E$.
En particulier, l'exemple \ref{E:Picard} ne peut pas se produire : les solutions
d'\'equations diff\'erentielles d'ordre $1$ ne rencontrent pas
de fronti\`ere naturelle au prolongement analytique.

\begin{remark}\rm
La projection d'une singularit\'e alg\'ebro\"\i de
de $f(x,x_0,y_0)$ d\'epend de la condition initiale $y_0$
et d\'ecrit tout $\Omega$ lorsque l'on d\'ecrit toutes les solutions.
Aussi, on les appelle {\it points critiques mobiles} dans
la litt\'erature.
\end{remark}

Introduisons maintenant les concepts g\'eom\'etriques
qui se cachent derri\`ere cet \'enonc\'e.
L'id\'ee de Paul Painlev\'e \'etait de consid\'erer l'\'equation diff\'erentielle $(E)$
comme un feuilletage singulier par courbes (\`a savoir les graphes des solutions)
dans le portrait de phases
$(x,y)\in\mathbb C\times\overline{\mathbb C}$ (Voir aussi [Re]).
Au voisinage de tout point $(x_0,y_0)\in\mathbb C^2$
satisfaisant $Q(x_0,y_0)\not=0$, les graphes des solutions locales $y=f(x)$
dont les conditions initiales $f(x_0)$ sont proches de $y_0$
sont les courbes de niveaux d'une submersion $H(x,y)$
(Th\'eor\`eme de Cauchy \`a param\`etre).
Pr\'ecis\'ement, il existe un voisinage
$U\ni(x_0,y_0)$ sur lequel $H:U\to\mathbb C$ est bien d\'efinie et r\'eguli\`ere
et telle que, pour toute condition initiale $(x_1,y_1)\in U$,
la solution locale
$f(x,x_1,y_1)$ est donn\'ee par l'\'equation fonctionnelle implicite
$H(x,f(x,x_1,y_1))=H(x_1,y_1)$. Ainsi, les graphes de solutions locales
d\'efinissent un feuilletage r\'egulier par courbes dans un voisinage de $(x_0,y_0)$.
Ce feuilletage $\mathcal F$ est globalement d\'efini sur le compl\'ementaire de $\{Q(x,y)=0\}$.
Il s'\'etend, le long de $\{Q(x,y)=0\}$, au voisinage de chaque point $(x_0,y_0)$ satisfaisant
$P(x_0,y_0)\not=0$ par les m\^emes arguments, apr\`es permutation du r\^ole des variables
$x$ et $y$, i.e. simplement en consid\'erant l'\'equation diff\'erentielle
${dx\over dy}={Q(x,y)\over P(x,y)}$.
Le feuilletage r\'egulier $\mathcal F$ est ainsi bien d\'efini sur le compl\'ementaire
$\mathbb C^2\setminus\{P(x,y)=Q(x,y)=0\}$. Quitte \`a diviser $P$ et $Q$ par un facteur commun,
nous pouvons toujours supposer ces polyn\^omes premiers entre eux.
Ainsi, l'ensemble singulier $\{P(x,y)=Q(x,y)=0\}$ consiste en un nombre fini de points,
\`a savoir les points d'ind\'etermination du second membre de $(E)$.
On dit que $\mathcal F$ est un feuilletage singulier \`a singularit\'es isol\'ees sur $\mathbb C^2$.

\eject

\begin{figure}[htbp]
\begin{center}

\input{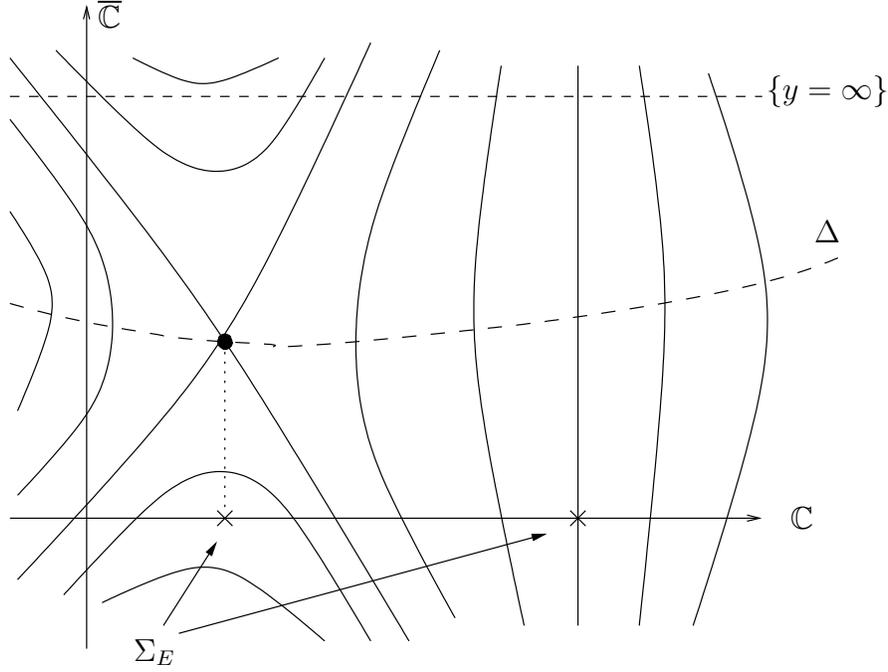}
 
\caption{Le feuilletage $\mathcal F$ et les singularit\'es fixes $\Sigma_E$}
\label{figure:7}
\end{center}
\end{figure}

Le feuilletage $\mathcal F$ s'\'etend en un feuilletage singulier \`a singularit\'es isol\'ees
sur $\mathbb C\times\overline{\mathbb C}$
de la mani\`ere suivante.
Supposons $P=P_m$ de degr\'e $m\in\mathbb N$ et $Q=Q_n$ de degr\'e $n\in\mathbb N$
en la variable $y$.
Apr\`es changement de variable $Y={1\over y}$, il vient :
$$(E)\ \ \ \ \ {dY\over dx}=-Y^2{P_m(x,{1\over Y})\over Q_n(x,{1\over Y})}={\widetilde P(x,Y)\over \widetilde Q(x,Y)}$$
o\`u $\widetilde P$ et $\widetilde Q$ sont les polynomes premiers entre eux donn\'es par :
$$\left\{\begin{matrix}
\widetilde P=-Y^{n+2}P_m(x,{1\over Y})&\text{et}&\widetilde Q=Y^nQ_n(x,{1\over Y})
&\text{si}& n+2\ge m\\
\widetilde P=-Y^mP_m(x,{1\over Y})&\text{et}&\widetilde Q=Y^{m-2}Q_n(x,{1\over Y})
&\text{si}&n+2<m
\end{matrix}\right.$$
On d\'eveloppe de nouveau les arguments pr\'ec\'edents pour \'etendre $\mathcal F$
en un feuilletage singulier \`a singularit\'es isol\'ees au voisinage de la droite \`a l'infini :
$$L_\infty=\mathbb C\times\overline{\mathbb C}\setminus\mathbb C^2=\{y=\infty\}.$$

\eject

On note $\Sigma_\mathcal F$ l'ensemble (fini) des singularit\'es de $\mathcal F$,
c'est \`a dire des points de $\mathbb C\times\overline{\mathbb C}$
au voisinage desquels $\mathcal F$ ne s'\'etend pas comme feuilletage r\'egulier :
$$\Sigma_\mathcal F\vert_{\mathbb C^2}=\{(x_0,y_0)\in\mathbb C^2\ ;\ P(x_0,y_0)=Q(x_0,y_0)=0\},$$
$$\Sigma_\mathcal F\vert_{L_\infty}=\{(x_0,\infty)\ ;\ \widetilde P(x_0,0)=\widetilde Q(x_0,0)=0\}.$$
On introduit aussi la courbe discriminante $\Delta$ comme l'ensemble des points de tangence
entre le feuilletage $\mathcal F$ et la fibration verticale :
$$\Delta\vert_{\mathbb C^2}=\{(x,y)\ ;\ Q(x,y)=0\}\ \ \ \text{et}\ \ \
\Delta\vert_{L_\infty}=\{(x,\infty)\ ;\ \widetilde Q(x,0)=0\}.$$

Nous allons distinguer les diff\'erentes positions possibles du feuilletage $\mathcal F$
vis \`a vis de la fibration verticale en un point
$(x_0,y_0)\in \mathbb C\times\overline{\mathbb C}$ suivant $5$ cat\'egories :
\begin{enumerate}
\item {\bf transversalit\'e} : $(x_0,y_0)\not\in\Delta$,
\item {\bf tangence simple} : $\mathcal F$ intersecte $\Delta$ transversalement en $(x_0,y_0)$,
\item {\bf tangence multiple} : $\mathcal F$ intersecte $\Delta$ avec multiplicit\'e en $(x_0,y_0)$,
\item {\bf singularit\'e} : $(x_0,y_0)\in\Sigma_\mathcal F$,
\item {\bf feuille verticale} : $Q(x_0,y)\equiv 0$,
i.e. $\{x=x_0\}$ est une feuille de $\mathcal F$.
\end{enumerate}
L'ensemble $\Sigma_E$ est la projection sur $\mathbb C$ des points de type (4) et (5) :
$$\Sigma_E=\left\{x_0\ ;\ Q(x_0,y)\equiv 0\ \ \text{ou}\ \ (x_0,y_0)\in\Sigma_\mathcal F\
\text{pour un}\ y_0\in\overline{\mathbb C}\ \right\}.$$

Dans les exemples \ref{E:puissanceComplexe} et \ref{E:SingNonIsolee}, 
$\Sigma_E$ est r\'eduit \`a $\{0\}$.
Dans l'exemple \ref{E:RiccatiRetournee}, on a $\Sigma_E=\{\zeta_1,\zeta_2,\zeta_3\}$
pour la premi\`ere \'equation
et toute solution $y=f(x)$ est une fonction multiforme
r\'eguli\`ere de $x$ sur $\Omega:=\mathbb C\setminus\{\zeta_1,\zeta_2,\zeta_3\}$~;
par contre, $\Sigma_E$ est vide pour la seconde \'equation
et toute solution $x=g(y)$
est une fonction alg\'ebro\"\i de transcendante de $x$
avec un ensemble dense de singularit\'es dans le plan.

\eject

\begin{theo}[II]
Soient $P,Q\in\mathbb C[x,y]$
et consid\'erons l'\'equation diff\'erentielle :
$$(E)\ \ \ \ \ {dy\over dx}={P(x,y)\over Q(x,y)}.$$
Notons $\Omega$
le compl\'ementaire dans $\mathbb C$ de l'ensemble $\Sigma_E$ donn\'e par le Th\'eor\`eme I.
Consid\'erons une solution analytique locale $f(x,x_0,y_0)$
et un chemin $\gamma:[0,1]\to\Omega$ issu de $x_0=\gamma(0)$.
Alors il existe une constante $\varepsilon>0$ telle que la fl\`eche
qui \`a une valeur initiale $y$ proche de $y_0$ fait correspondre
l'ensemble fini de valeurs :
$$\varphi_\gamma(y):=\left\{f_{\gamma'}(x_1,x_0,y)\ ;\ \begin{matrix}
\gamma':[0,1]\to\Omega\ \text{est $\varepsilon$-proche de $\gamma$ et}\\
\text{\'evite les singularit\'es interm\'ediaires}
\end{matrix}\right\}$$
d\'efinit, sur tout disque $D$ suffisamment petit centr\'e en $y_0$,
une fonction multivalu\'ee de $y$
donn\'ee par un nombre fini de fonctions alg\'ebro\"\i des (finies) sur $D$.
\end{theo}

Si le concept de feuilletage n'est pas clairement d\'efini dans [Pa],
le Th\'eor\`eme II est bien connu des ``feuilleteurs'' dans le contexte suivant.

\begin{remark}\rm
Le choix d'un prolongement alg\'ebro\"\i de de $f(x,x_0,y_0)$ correspond au choix
d'un chemin $\widetilde\gamma:[0,1]\to \Omega\times\overline{\mathbb C}$ relevant $\gamma$
dans la feuille $\mathcal F_{p_0}$ passant par le point $p_0=(x_0,y_0)$ :
$\widetilde\gamma(0)=p_0$ et $\pi\circ\widetilde\gamma=\gamma$ (ici
$\pi$ est la projection verticale). Il appara\^\i tra clairement
dans la suite que le prolongement obtenu est analytique
(i.e. non singulier en $x_1$) d\`es que
le feuilletage $\mathcal F$ est transverse \`a la verticale $T_1:=\{x=x_1\}$
au point d'arriv\'ee $\widetilde\gamma(1)$. En d'autres termes,
{\it \'etant donn\'e un chemin $\widetilde\gamma:[0,1]\to\Omega\times\overline{\mathbb C}$
\`a valeurs dans une feuille joignant deux points
$\gamma(i)$, $i=0,1$, en lesquels le feuilletage est transverse \`a la verticale $T_i$,
il existe un $\varepsilon>0$
tel que la fl\`eche $\widetilde\gamma'(0)\mapsto\widetilde\gamma'(1)$
qui fait se correspondre les extr\'emit\'es de tous les chemins $\widetilde\gamma'$
tangents au feuilletage, joignant $T_0$ \`a $T_1$ et $\varepsilon$-proches
de $\widetilde\gamma$
d\'efinit une application uniforme et holomorphe sur un voisinage de $y_0$.}
C'est le lemme fondamental de l'holonomie.
\end{remark}

Comme l'a maintes fois rappel\'e Painlev\'e (voir [Bo], p.142),
il est important de comprendre
que l'\'enonc\'e pr\'ec\'edent n'est que local dans la variable $y$.
Le prolongement analytique de la fl\`eche $y\mapsto f(x_1,x_0,y)$
peut conduire \`a des singularit\'es non alg\'ebro\"\i des dans le plan des $y$
(voir \S \ref{S:ExemplePainleve}).

\eject

\section{Preuve du Th\'eor\`eme I }

La d\'emonstration du premier th\'eor\`eme repose sur la compacit\'e
de la fibre verticale $\overline{\mathbb C}$.
Nous commen\c cons par prouver un cas particulier important.

\begin{prop}\label{P:Riccati}
Soient $P,Q\in\mathbb C[x,y]$
et consid\'erons l'\'equation diff\'erentielle :
$$(E)\ \ \ \ \ {dy\over dx}={P(x,y)\over Q(x,y)}.$$
S'il existe un point $x_0\in\mathbb C$ au dessus duquel les solutions
$y(x)$ de l'\'equation $(E)$ n'ont pas d'autre singularit\'e
que des p\^oles, alors $(E)$ est une \'equation de Riccati :
$${dy\over dx}={a(x)y^2+b(x)y+c(x)\over Q(x)}$$
o\`u $a,b,c,Q\in\mathbb C[x]$ sont des polyn\^omes d'une variable.
De plus, on a $\Sigma_E=\{Q(x)=0\}$ et toute solution $y(x)$ est
m\'eromorphe sur le rev\^etement universel 
$\widetilde{\mathbb C\setminus\Sigma_E}$.
\end{prop}

\begin{proof}L'absence de singularit\'e alg\'ebro\"\i de
autre que des p\^oles au dessus de $x_0$ se traduit
par la transversalit\'e de la droite verticale $\{x=x_0\}$ 
avec le feuilletage $\mathcal F$ dans $\mathbb C\times\overline{\mathbb C}$.
En particulier, $Q(x_0,y)\not=0$ pour tout $y\in\mathbb C$,
c'est \`a dire $Q=q_0(x)+q_1(x)y+\cdots+q_n(x)y^n$
avec $q_0(x_0)\not=0$ et $q_k(x_0)=0$ pour $k=1,\ldots,n$. Apr\`es changement de coordonn\'ees
$Y:=1/y$, il appara\^\i t que $\widetilde Q(x_0,Y)\equiv q_0(x_0)Y^k$ o\`u $k:=\sup\{m-2,n\}$.
La transversalit\'e de $\{x=x_0\}$ \`a $y_0=\infty$ se traduit par $\widetilde Q(x_0,0)\not=0$.
Par suite, $n=0$ et $m\le 2$ ce qui nous conduit \`a l'\'equation de Riccati. Dans ce cas, toute fibre verticale $\{x=x_0\}$, $Q(x_0)\not=0$, est transverse au feuilletage.
Par cons\'equent, la premi\`ere projection $(x,y)\mapsto x$
induit sur chaque feuille $L$ un rev\^etement 
$L\to\mathbb C\setminus\Sigma_E$. Ainsi, $L$ est le graphe d'une
fonction n'ayant que des p\^oles au dessus de 
$\mathbb C\setminus\Sigma_E$, \`a savoir les points d'intersection
de $L$ avec la droite \`a l'infini $\{y=\infty\}$.
\end{proof}

\eject

\begin{remark}\rm
Une \'equation diff\'erentielle dont les solutions sont toutes m\'ero\-morphes
est une \'equation de Riccati.
En effet, les autres \'equations poss\`edent toujours des solutions avec singularit\'es
alg\'ebro\"\i des multiformes (points d'intersection du graphe
avec $\Delta$).
Cependant, tr\`es peu d'\'equations de Riccati ont des solutions m\'eromorphes.
L'exemple \ref{E:RiccatiRetournee}~:
$${dy\over dx}={\alpha_1\over x-\zeta_1}+{\alpha_2\over x-\zeta_2}+{\alpha_3\over x-\zeta_3}$$
est une \'equation de Riccati dont toute solution est multiforme.
\end{remark}

La proposition \ref{P:Riccati} repose sur la trivialisation locale des points de type (1).
Pour les points $(x_0,y_0)$ de type (2),
on peut encore donner des mod\`eles locaux :

\begin{lem}\label{L:tangenceSimple}
Soit $\mathcal F$ un feuilletage r\'egulier d\'efini au voisinage
de l'origine $\underline 0\in\mathbb C^2$ par une submersion $H:U\to\mathbb C$,
avec une tangence verticale
en $\underline 0$. Supposons que le discriminant $\Delta:=\{{\partial H\over \partial y}=0\}$
soit transverse \`a la verticale $\{x=0\}$ en $\underline 0$.
Alors, il existe un changement de coordonn\'ees locales $\Phi(x,y)=(x,\phi(x,y))$ pr\'eservant
la  fibration verticale et envoyant $\mathcal F$ sur le feuilletage d\'efini par $H_k(x,y)=x+y^k$,
o\`u $(k-1)\in\mathbb N^*$ est l'ordre de ${\partial H\over \partial y}$ le long de $\Delta$.
En d'autres termes, il existe un diff\'eomorphisme local $\varphi$ tel que
$H=\varphi\circ H_k\circ\Phi$.
\end{lem}

\begin{proof}Par un changement de coordonn\'ees pr\'eliminaire (pr\'eservant la fibration verticale),
on peut d\'ej\`a supposer $\Delta:\{y=0\}$. Alors, ${\partial H\over \partial y}=y^{k-1} u(x,y)$
pour un entier $(k-1)\in\mathbb N^*$ et une unit\'e holomorphe $u(x,y)$, $u(0,0)\not=0$.
Par int\'egration, on obtient $H=f(x)+y^k\widetilde u(x,y)$ avec une nouvelle unit\'e
$\widetilde u(x,y)$
et une fonction $f(x)$ que l'on peut supposer s'annuler en $0$, i.e. $H(0,0)=0$.
Puisque $H$ est une submersion transverse \`a $\Delta$, $\varphi(x)=H(x,0)$
est encore une submersion sur $\Delta$. Donc, $\varphi$ est un diff\'eomorphisme local
et $\widetilde H:=\varphi^{-1}\circ H$ est une nouvelle submersion d\'efinissant $\mathcal F$
de la forme $\widetilde H(x,y)=x+y^k\widetilde u(x,y)$.
Finalement, si $v(x,y)$ d\'esigne une racine primitive $k^{\text{i\`eme}}$ de $\widetilde u(x,y)$,
alors le changement de coordonn\'ees d\'efini par $\phi(x,y)=y\cdot v(x,y)$ satisfait
$\widetilde H(x,y)=x+(\phi(x,y))^k=H_k\circ \Phi(x,y)$.
\end{proof}

\eject

\begin{prop}
Soit $T:=\{x=x_0\}\subset\mathbb C\times\overline{\mathbb C}$ une verticale
transverse \`a $\Delta$ (tout point le long de $T$ est de type (1) ou (2)). Alors, on peut
trouver un disque $D\subset\mathbb C$ centr\'e en $x_0$ tel que la projection verticale
restreinte au voisinage tubulaire $D\times\overline{\mathbb C}$ induise un rev\^etement ramifi\'e
des feuilles (apr\`es restriction) sur $D$. De plus, chacun de ces rev\^etements ramifi\'es
a au plus un branchement. Les multiplicit\'es de branchement possibles sont les multiplicit\'es
$k_1,\ldots,k_n\in\mathbb N^*$ de
$Q(x,y)$ et $\widetilde Q(x,Y)$
le long des composantes irr\'eductibles non verticales de $\Delta$.
En d'autres termes, toute solution $f(x,x_0,y_0)$ dans $x\in D$ est ou bien uniforme,
ou bien alg\'ebro\"\i de avec exactement une singularit\'e et $k_i$ d\'eterminations
pour un $i\in\{1,\ldots,n\}$.
\end{prop}

\begin{proof}Elle se d\'eduit par un argument de compacit\'e appliqu\'e aux mod\`eles
locaux donn\'es le long de la verticale par le Lemme \ref{L:tangenceSimple}.
\end{proof}

\begin{example}\rm
Le feuilletage globalement donn\'e par $H_k=x+y^k$ correspond
\`a l'\'equation diff\'erentielle ${dy\over dx}=-{1\over ky^{k-1}}$.
Toute solution $y=(x_0+y_0^k-x)^{1/k}$ poss\`ede exactement une singularit\'e alg\'ebro\"\i de
en $x_1:=x_0+y_0^k$. La position de cette singularit\'e d\'epend des conditions initiales~:
$x_1$ est un {\it point critique mobile}.
\end{example}

\begin{example}\rm
Pour $n>2$, l'\'equation diff\'erentielle :
$$(E)\ \ \ \ \ {dy\over dx}={a_n(x)y^n+\cdots+a_1(x)y+a_0(x)\over Q(x)},$$
o\`u $a_i$ et $Q$ d\'esignent des polyn\^omes en $x$,
a pour discriminant $\Delta=L_\infty\cup\{Q(x)=0\}$ avec multiplicit\'e $k=n-1$
le long de $L_\infty$.
\end{example}

\begin{example}\label{E:tangenceDegeneree}\rm
Le feuilletage globalement d\'efini par $H={x-y^3\over 1-y}$ correspond
\`a l'\'equation diff\'erentielle ${dy\over dx}={y-1\over 2y^3-3y^2+x}$. Son discriminant
poss\`ede une tangente verticale \`a l'origine. La solution $f(x,0,0)=x^{1/3}$ poss\`ede
exactement une singularit\'e alg\'ebro\"\i de en $0$ avec $3$ d\'eterminations.
Puisque $y=f(x,x_0,y_0)$ est l'inverse de $x=g(y)=y^3-cy+c$, $c=H(x_0,y_0)$,
elle a $2$ singularit\'es alg\'ebro\"\i des avec $2$ d\'eterminations locales
autour de chacune d'elles ; d'un point de vue global, elle poss\`ede $3$ d\'eterminations.
Lorsque $c\to 0$, ces deux
singularit\'es bifurquent en une seule d'ordre plus grand.
\end{example}

\begin{remark}\rm
D\`es que l'on autorise la pr\'esence de singularit\'es dans $\Delta$,
on ne peut plus esp\'erer obtenir de mod\`eles polynomiaux.
En effet, si l'on consid\`ere le cas le plus simple o\`u $\Delta$ est \`a croisements ordinaires,
disons $\Delta=\{y^2-x^2=0\}$,
alors on peut mettre en \'evidence l'invariant suivant. Notons $\Delta_1=\{y-x=0\}$
et $\Delta_2=\{y+x=0\}$ et d\'efinissons un germe de diff\'eomorphisme
$h:\Delta_1\to\Delta_1$ de la fa\c con suivante. Envoyons d'abord $p\in\Delta_1$
sur un point $p'\in\Delta_2$  par projection verticale, c'est \`a dire suivant le feuilletage vertical,
puis envoyons $p'$ sur l'unique point
$h(p)\in\Delta_1$ qui est dans la m\^eme feuille pour $\mathcal F$.
Si $H$ est une submersion d\'efinissant le feuilletage $\mathcal F$, alors $h(x)=h_1^{-1}\circ h_2$
o\`u $h_1(x)=H(x,x)$ et $h_2(x)=H(x,-x)$. Apr\`es changement de coordonn\'ees pr\'eservant
le feuilletage vertical, c'est \`a dire de la forme $\Phi(x,y)=(\varphi(x),\phi(x,y))$,
$h$ se trouve conjugu\'e par $\varphi$.
Par exemple, la submersion $H=x+({y^3\over 3}-x^2y+{2\over 3}x^3)u(x)$,
$u(0)=0$, est de discriminant $\Delta=\{y^2-x^2=0\}$ et r\'ealise le diff\'eomorphisme
$h(x)=x+{4\over 3}x^3u(x)$. Donc, tout diff\'eomorphisme $h$ tangent
\`a l'identit\'e \`a l'ordre $2$ est r\'ealisable comme invariant local d'un tel feuilletage.
La classification de tels diff\'eomorphismes modulo conjugaison analytique donne naissance
aux modules d'\'Ecalle-Malgrange-Voronin qui sont fonctionnels,
plus grand que l'espace des polyn\^omes.
\end{remark}

\begin{lem}\label{L:tangenceGenerale}
Soit $\mathcal F$ un feuilletage r\'egulier d\'efini au
voisinage de $(0,0)\in\mathbb C^2$ par une submersion $H:U\to\mathbb C$, $H(0,0)=0$.
Supposons que $\{H(x,y)=0\}$ ne soit pas la droite verticale $\{x=0\}$ :
${\partial H\over\partial y}(0,y)$ est d'ordre fini $k\in\mathbb N$ en $y=0$.
Alors il existe un disque $D$ centr\'e en $0\in\mathbb C$ et $\varepsilon>0$
tels que, pour $\vert c\vert\le\varepsilon$, les solutions $y=f(x)$
de l'\'equation implicite $H(x,y)=c$
sont alg\'ebro\"\i des sur $D$ avec $\le k+1$ d\'eterminations et au plus $k$ singularit\'es.
Plus pr\'ecis\'ement, pour $c=0$, il y a seulement une solution $y(x)$
avec seulement une singularit\'e en $x=0$
et $k+1$ d\'eterminations autour. Pour $c\not=0$ suffisament proche de $0$,
il y a encore une unique solution $y(x)$ avec $\le k$ singularit\'es $\not=0$ et
$k+1$ d\'eterminations autour.
Alors que $\vert c\vert$ cro\^\i t, les singularit\'es s'\'echappent par le bord de $D$
faisant cro\^\i tre le nombre de solutions distinctes. Finalement,
pour $\vert c\vert=\varepsilon$, il y a exactement $k+1$ solutions distinctes uniformes.
\end{lem}

\begin{proof}Lorsque $k=0$, le th\'eor\`eme des fonctions implicites
donne des solutions uniformes $\varphi_c(x)$
analytiques sur un polydisque $(x,c)\in D\times D_\varepsilon$ centr\'e en $(0,0)$.
\`A partir de maintenant, supposons $k\ge 1$.
Le th\'eor\`eme des fonctions implicites
donne des solutions uniformes $x=\psi_c(y)$
\`a l'\'equation fonctionnelle $H(\psi_c(y),y)=c$ analytiques sur un polydisque
$(y,c)\in D'\times D_\varepsilon$ centr\'e en $(0,0)$. Pour $c=0$, $\psi_0(y)$ poss\`ede
un point critique s'annulant \`a l'ordre $k+1$ en $y=0$. Quitte \`a diminuer $D'$,
$\psi_0(y)$ est conjugu\'e,
sur un voisinage de l'adh\'erence $\overline{D'}$, \`a la fonction $y\mapsto y^{k+1}$.
En particulier, il existe $2$ disques $D^-\subset D^+$ centr\'es en $x=0$ tels que
la fronti\`ere $\partial D'$ soit envoy\'ee par $\psi_0$ dans l'anneau $D^+\setminus \overline{D^-}$
avec indice $k+1$. Quitte \`a diminuer $\varepsilon$, $\psi_c$ envoit encore $\partial D'$ dans
$D^+\setminus \overline{D^-}$ avec indice $k+1$ et,  par le th\'eor\`eme de Rouch\'e,
poss\`ede encore $k$
points critiques (compt\'es avec multiplicit\'e) dans $D'$ pour $\vert c\vert\le\varepsilon$.
Ceci signifie que pour pour tout sousdisque $D\subset D^-$, $\psi_c$ induit un rev\^etement ramifi\'e
$\psi_c^{-1}(D)\to D$ de degr\'e $k+1$ avec au plus $k$ points critiques.
Les solutions alg\'ebro\"\i des
$y(x)$ de l'\'enonc\'e sont les applications inverses de $\psi_c$, allant de $D$ vers les
composantes connexes de $\psi_c^{-1}(D)$ ce qui prouve le premier point.

Il reste \`a montrer que les images des points critiques de $\psi_c$, qui deviennent
les points singuliers des solutions $y(x)$, s'\'echappent de $0$ et, par suite,
de tout disque $D$ suffisamment petit
lorsque $\vert c\vert$ cro\^\i t. En fait, nous allons montrer que les singularit\'es
des solutions $y(x)$ sont des germes de fonctions alg\'ebro\"\i des non constantes de $c$ en $0$.
Remarquons d'abord que les points critiques $y_1(c),\ldots,y_k(c)$ de $\psi_c$ sont
les d\'eterminations de fonctions alg\'ebro\"\i des de $c$ en $c=0$.

En effet, notons $\Delta$ la courbe discriminante de $\mathcal F$,
$\Delta=\{{\partial H\over\partial y}=0\}$.
Dans les nouvelles coordonn\'ees $(c,y):=(H(x,y),y)$, la fonction $H$ devient
la premi\`ere coordonn\'ee et
l'ancienne coordonn\'ee $x=x(c,y)$ devient une submersion de discriminant $\Delta$.
Alors, les points critiques $y_1(c),\ldots,y_k(c)$ sont les points d'intersection
de la droite verticale $\{H=c\}$
avec la courbe $\Delta$. Donc, les points critiques sont les param\'etrisations de Puiseux
de $\Delta$
dans ces coordonn\'ees et sont par suite alg\'ebro\"\i des au voisinage de $0$.
Finalement, les points singuliers $\zeta_i(c)$ des solutions $y(x)$ sont donn\'es par
$\zeta_i(c)=\psi_c(y_i(c))$ et sont aussi alg\'ebro\"\i des (composition de
d'une s\'erie de Puiseux par une s\'erie enti\`ere).

Puisque $H$ poss\`ede une tangente isol\'ee
avec la droite verticale $\{x=0\}$ en $(0,0)$, le graphe $\{H(x,y)=c\}$ de $\psi_c$
intersecte transversalement la droite $\{x=0\}$ pour $c\not=0$ et tout fonction inverse
correspondante $y(x)$ est analytique en $0$. Donc, aucun de ces germes alg\'ebro\"\i des $\zeta_i(c)$
ne peut \^etre constant $\equiv0$. En particulier, pour $\vert c\vert=\varepsilon$,
chacune de ces d\'eterminations s'est \'echapp\'e de tout disque suffisamment petit $D$.
\end{proof}

\begin{remark}\rm Dans l'exemple \ref{E:tangenceDegeneree}, 
les points critiques des solutions $x(y)$ sont
les $2$ d\'eterminations de la fonction alg\'ebro\"\i de ${c^{1/2}\over\sqrt{3}}=\{y_1(c),y_2(c)\}$.
Les points singuliers de l'inverse $y(x)$ sont les valeurs critiques de la fonction pr\'ec\'edente,
c'est \`a dire sont les $2$ d\'eterminations de la fonction alg\'ebro\"\i de
$c-{2\over3\sqrt{3}}c^{3/2}=\{\zeta_1(c),\zeta_2(c)\}$.
Dans la remarque 5.8, les points singuliers sont $2$ fonctions uniformes donn\'ees
par $\zeta_1(c)=h_1^{-1}(x)$ et $\zeta_2(c)=h_2^{-1}(x)$.
L'invariant $h$ discut\'e dans la remarque 5.8 n'est autre que
l'application permuttant les singularit\'es. Aussi, dans l'exemple \ref{E:tangenceDegeneree},
si $\Delta$ est param\'etr\'e  par $z\mapsto(3z^2-2z^3,z)$, alors $H\vert_{\Delta}(z)=z^2$.
Donc, l'invariant $h$, vu comme application $\Delta\to\Delta$ est simplement $z\mapsto-z$
alors qu'il devient alg\'ebro\"\i de lorsqu'il est vu dans la variable $x$.
\end{remark}

Nous d\'eduisons finalement du Lemme \ref{L:tangenceGenerale} la version
quantitative suivante du Th\'eor\`eme I.

\begin{cor}
\'Etant donn\'es l'\'equation diff\'erentielle :
$$(E)\ \ \ \ \ {dy\over dx}={P(x,y)\over Q(x,y)}$$
et $\Sigma_E\in\mathbb C$ d\'efinis comme pr\'ec\'edemment,
il existe un $k\in\mathbb N$ tel que toute solution est alg\'ebro\"\i de
sur $\Omega=\mathbb C\setminus\Sigma_E$ avec singularit\'es ramif\'ees \`a l'ordre $\le k$.
De plus, pour tout chemin diff\'erentiable $\gamma:[0,1]\to\Omega$ et
toute solution analytique locale $f(x,x_0,y_0)$, $x_0=\gamma(0)$, le nombre
de d\'eterminations obtenues en $x_1=\gamma(1)$ par prolongement analytique le long de $\gamma$
contournant les singularit\'es interm\'ediaires est born\'e par
$(k+1)^{\vert\gamma\vert/\varepsilon(\gamma)}$
o\`u $\vert\gamma\vert$ d\'esigne la longueur de $\gamma$ et $\varepsilon(\gamma)>0$ est
une fonction d\'ecroissante de la distance $d(\gamma([0,1]),\Sigma_E)$ de $\gamma$
\`a l'ensemble singulier $\Sigma_E$
(pour la distance sph\'erique sur $\overline{\mathbb C}\supset\mathbb C$).
\end{cor}

\eject

\begin{proof}
Pr\`es de tout point $(x_0,y_0)\in\mathbb C^2$ satisfaisant $x_0\not\in\Sigma_E$,
le lemme 7 fournit une description des solutions de $(E)$ sur un voisinage ouvert
de la forme $\{\vert x\vert<\varepsilon_1,\vert H\vert<\varepsilon_2\}$,
$\varepsilon_1,\varepsilon_2>0$.
Par un argument de compacit\'e, la fibre $\{x=x_0\}$
peut \^etre recouverte par un nombre fini de tels voisinages. Il existe donc un disque $D(x_0)$
centr\'e en $x_0$ tel que toute solution $y(x)$ de $(E)$ dans $D$ est alg\'ebro\"\i de
avec au plus $k$ points singuliers et au plus $k+1$ d\'eterminations. L'entier $k$ est
la multiplicit\'e totale d'intersection entre $\Delta$ et la droite verticale $\{x=x_0\}$
et est alors ind\'ependant de $x_0$. Tout sous-ensemble compact $K\subset\Omega$
(par exemple $K=K_r:=\{x\ ;\ d(x,\Sigma_E)\ge r\}$, $r>0$) est recouvert par un nombre fini
de tels disques $D_1,\ldots,D_N$. Alors pour tout chemin $\gamma:[0,1]\to K$ et pour toute
solution analytique locale $f(x,x_0,y_0)$, $x_0=\gamma(0)$, on peut effectuer le prolongement
alg\'ebro\"\i de le long de $\gamma$ en recollant ensemble un nombre fini de solutions
alg\'ebro\"\i des dans les disques $D_i$ rencontr\'e. Le nombre maximal de disques successifs
$D(x_i)$ n\'ecessaire pour recouvrir $\gamma$ est born\'e par la longueur $\vert\gamma\vert$.
En effet, fixons un $\varepsilon>0$ tel que les disques $\varepsilon$-rogn\'es
$D_i^-:=\{x\in D_i\ ;\ d(x,\partial D_i)>\varepsilon\}$
recouvrent encore $K$. Alors si $\gamma(t_0)$ est dans un $D_i^-$, $\gamma(t)$ reste dans $D_i$
pour $t\in[t_0,t_0+\varepsilon]$ avant de changer pour un autre disque $D_j^-$.
Ainsi, on a besoin de changer au plus ${\vert\gamma\vert\over\varepsilon}$ fois de disques $D_i$
pour pouvoir effectuer le prolongement alg\'ebro\"\i de le long de $\gamma$.
En d'autres termes, au plus
${\vert\gamma\vert\over\varepsilon}+1$ disques distincts ont \'et\'e rencontr\'es
durant le prolongement analytique et au plus $(k+1)^{({\vert\gamma\vert\over\varepsilon}+1)}$
d\'eterminations distinctes peuvent \^etre obtenues pour $f(x_1,x_0,y_0)$, $x_1=\gamma(1)$.
\end{proof}

\eject

\section{Preuve du th\'eor\`eme II et l'exemple de Painlev\'e}\label{S:ExemplePainleve}

\begin{proof}[Preuve du th\'eor\`eme II]
Soient $f(x,x_0,y_0)$ et $\gamma:[0,1]\to \Omega$ comme dans l'\'enonc\'e
et soit $\Sigma_\gamma(f)$ l'ensemble des singularit\'es interm\'ediaires possibles
(voir propri\'et\'e (B) du \S 3).
Tout prolongement alg\'ebro\"\i de $f_\gamma$ va donner naissance \`a une composante irr\'eductible
de $\varphi$. Un tel prolongement est caract\'eris\'e, d'apr\`es la remarque 14,
par un rel\`evement :
$$\widetilde\gamma:[0,1]\to L_0\subset\mathbb C\times\overline{\mathbb C},\
\Pi_1\circ\widetilde\gamma=\gamma\ \ \ \text{avec}\ \ \ \widetilde\gamma(0)=(x_0,y_0)$$
dans la feuille $L_0$ passant par $(x_0,y_0)$.
Soient $\widetilde\gamma_1,\ldots,\widetilde\gamma_N:[0,1]\to\Omega\times\overline{\mathbb C}$
tous les rel\`evements possibles de $\gamma$, caract\'erisant tous les prolongements
alg\'ebro\"\i des de $f$.
En particulier, tous ces chemins co\"\i ncident pour $t$ petit et bifurquent
au f\^ur et \`a mesure que l'on rencontre les points de $\Sigma_\gamma(f)$.

D'apr\`es le lemme 7 et par compacit\'e de la fibre, pour tout $t\in[0,1]$, on peut trouver un
disque $D$ centr\'e en $\gamma(t)$ tel que toute feuille $L$ du feuilletage $\mathcal F\vert_D$
restreinte au cylindre $D\times\overline{\mathbb C}$ est un rev\^etement ramifi\'e sur $D$.
Notons $L_l$ la feuille locale contenant $\widetilde\gamma_l(t)$. Quitte \`a diminuer
$D$, on peut supposer que chaque application $\pi:L_l\to D$ est ou bien un diff\'eomorphisme,
ou bien un rev\^etement ramifi\'e ne ramifiant qu'au dessus de $\widetilde\gamma_l(t)$.
Dans ce dernier cas, remarquons que
$\gamma(t)\in\Sigma_\gamma(f)$ et, si l'ordre de branchement est $q$, alors $q$ rel\`evements
$\widetilde\gamma_l$ coincident jusqu'\`a $t$ et deviennent distincts juste apr\`es.
Par compacit\'e du chemin $\gamma$, on peut extraire de ces disques un recouvrement fini,
disons :
$$D_k\supset\gamma([t_k,t_{k+1}])\ \ \ \text{o\`u}\ \ \ 0=t_0<t_1<\cdots<t_n<t_{n+1}=1.$$

Pour $k=0,\ldots,n$ et $l=1,\ldots,N$, notons $L_{k,l}$ la feuille du
feuilletage restreint $\mathcal F\vert_{D_k}$ contenant $\widetilde\gamma_l(t_k)$.
En particulier, $L_{k,l}\cap L_{k+1,l}$ est envoy\'e diff\'eomorphique\-ment sur $D_k\cap D_{k+1}$
par $\pi$.
Soit $H_{k,l}:U_{k,l}\to\mathbb C$ une submersion definissant le feuilletage sur un
voisinage tubulaire $U_{k,l}$ de $L_{k,l}$. Apr\`es composition (\`a gauche) des $H_{k,l}$
par des diff\'eomorphismes convenables, on peut supposer $H_{0,l}(x_0,y)=y$ et
$H_{k+1,l}\equiv H_{k,l}$
au voisinage de $\widetilde\gamma_l(t_{k+1})$. Pour $\varepsilon>0$ suffisamment petit,
on peut supposer que, pour tout $\vert y_0'-y_0\vert<\varepsilon$, l'ensemble
$\{x\in D_k\ \text{et}\ H_{k,l}(x,y)=y_0'\}$ n'est pas vide et est une feuille compl\`ete pour
le feuilletage restreint $\mathcal F\vert_{D_k}$. Pour simplifier, rempla\c cons $U_{k,l}$
par l'ouvert
$\{x\in D_k\ \text{et}\ \vert H_{k,l}(x,y)-y_0\vert<\varepsilon\}$.

\eject

Soit $\gamma':[0,1]\to \Omega$ une $\varepsilon$-perturbation de $\gamma$
avec $\gamma'(0)=x_0$. Si $\varepsilon>0$ est suffisamment petit,
Alors $D_k$ reste un recouvrement ordonn\'e pour $\gamma'$. De plus,
on peut supposer $\gamma'(t_k)\in D_k\cap D_{k+1}$.
Soit $\widetilde\gamma':[0,1]\to \mathbb C\times\overline{\mathbb C}$ un rel\`evement
de $\gamma'$ dans le feuilletage correspondant \`a un prolongement alg\'ebro\"\i de
d'une solution $f(x,x_0,y_0')$ avec $\vert y_0'-y_0\vert<\varepsilon$.
Il s'en suit que $\widetilde\gamma'(1)=(x_1,y_1')$ o\`u $x_1=\gamma(1)$
et $y_1'$ est l'une des valeurs prises par $\varphi_\gamma$ en $y_0'$.
De plus, toute valeur de $\varphi_\gamma$ est obtenue de cette mani\`ere.

\`A pr\'esent, nous sommes \`a m\^eme de prouver que $\widetilde\gamma'$
reste proche de $\widetilde\gamma$.
En effet, supposons que $\widetilde\gamma'(t_k)$ soit contenu dans un $U_{k,l}$.
Alors $\widetilde\gamma'(t_k)$ est dans la feuille locale
$L_k:=\{x\in D_k\ \text{et}\ H_{k,l}(x,y)=y_0'\}$.
Puisque $\gamma'([t_k,t_{k+1}])$ est contenu dans $D_k$, $\widetilde\gamma'([t_k,t_{k+1}])$
est compl\`etement contenu dans $L_k$. En particulier, quitte \`a remplacer l'indice $l$
par un autre convenable parmis ceux pour lesquels
$H_{k,l}:U_{k,l}\to\mathbb C$ reste inchang\'e,
$\widetilde\gamma'(t_{k+1})$ est dans $U_{k,l}\cap U_{k+1,l}$.
Par r\'ecurrence, puisque $\widetilde\gamma'(0)\in U_{0,1}$, il s'en suit
qu'il existe un $l=1,\ldots,N$ tel que $\widetilde\gamma'(t_k)\in U_{k,l}$ pour tout
$k=0,\ldots,n$. Finalement,
$\widetilde\gamma'(1)=(x_1,y_1')\in L_n=\{x\in D_n\ \text{et}\ H_{k,l}(x,y)=y_0'\}$.
Par cons\'equent, $y_1'$ est la valeur prise en $y_0'$ par une des d\'eterminations
des fonctions alg\'ebro\"\i des
$\varphi_l(y)$ d\'efinies pour $\vert y-y_0\vert<\varepsilon$  par
$H_{n,l}(x_1,\varphi_l(y))=y$.\end{proof}

L'\'enonc\'e original du th\'eor\`eme II parait ambig\"u aujourd'hui.
Dans les  {\it Le\c cons de Stokholm}, il est \'ecrit (voir [Pa], tome 1, p.210) :

\noindent {\it Soit $y_0$ la valeur de $y(x)$ pour $x=x_0$, et soit
$y=\varphi(x,y_0,x_0)$ l'int\'egrale g\'en\'erale de (A). Si $\overline{x}$, $\overline{x_0}$
d\'esignent deux valeurs num\'eriques quelconques distinctes des valeurs $\xi$,
la fonction $y=\varphi(\overline{x},y_0,\overline{x_0})$ ne pr\'esente dans tout le plan
des $y_0$ (\`a distance finie ou infinie) que des points singuliers alg\'ebriques.}

On aurait tendance \`a penser, en lisant ce texte, que pour tous points $x_0$ et $x_1$
distincts de $\Sigma_E$, tout chemin $\gamma$ les joignant et toute valeur
$y_1=f(x_1,x_0,y_0)$ obtenue en choisissant un relev\'e de $\gamma$ dans le feuilletage,
le germe de fonction $\varphi:(\mathbb C,y_0)\to(\mathbb C,y_1)$ correspondant d\'efini par
la fl\`eche $y\mapsto f(x_1,x_0,y)$ admet un prolongement alg\'ebro\"\i de
sur $\overline{\mathbb C}$. Cependant, la simple lecture de la preuve originale de Painlev\'e
nous montre qu'il n'en est rien : le th\'eor\`eme prouv\'e est purement local.
De plus, il est tr\`es facile de se convaincre que le germe $\varphi$ poss\`ede
en g\'en\'eral des singularit\'es non alg\'ebro\"\i des.

\eject

\begin{example}\label{E:Dulac}\rm
 Si $P$ et $Q$
sont deux polyn\^omes de degr\'e $d\in\mathbb N^*$ suffisamment g\'en\'eriques,
alors le feuilletage induit par l'\'equation diff\'erentielle correspondante :
$$(E)\ \ \ \ \ {dy\over dx}={P(x,y)\over Q(x,y)}$$
poss\`ede exactement $d^2$ points singuliers dans $\mathbb C^2$,
tous hyperboliques.
D'apr\`es le th\'eor\`eme de Poincar\'e, il existe, au voisinage d'un tel point $p\in\mathbb C^2$,
un syst\`eme de cordonn\'ees analytiques locales redressant le feuilletage
sur celui d\'efini par :
$${dy\over dx}=\alpha{y\over x}\ \ \ \text{c'est \`a dire}\ \ \ \{y\cdot x^{-\alpha}=constante\}$$
pour un nombre complexe $\alpha\not\in\mathbb R$.
En particulier, le feuilletage poss\`ede deux courbes locales invariantes,
lisses et transverses, au point $p$ qui n'ont aucune raison, ni l'une ni l'autre, d'\^etre verticales.
Ces deux feuille vont imanquablement intersecter n'importe quelle droite verticale $T=\{x=x_0\}$.
Notons $L_0$ et $L_1$ ces deux feuilles, choisissons (le cas \'ech\'eant)
un point d'inter\-sec\-tion respectif $p_0=(x_0,y_0)\in L_0\cap T$ et $p_1=(x_0,y_1)\in L_1\cap T$
ainsi qu'un chemin $\widetilde\gamma_i$ joignant $p_i$ \`a $p$ dans $L_i$ pour $i=1,2$.
La concat\'enation $\widetilde\gamma:=\widetilde\gamma_0\cdot \widetilde\gamma_1^{-1}$ est un chemin
allant de $p_0$ \`a $p_1$, tangent au feuilletage et changeant de feuille au point $p$.
Au vu du mod\`ele local, le germe d'application retour $\varphi:(\mathbb C,y_0)\to(\mathbb C,y_1)$
correspondant (d\'efini comme dans l'\'enonc\'e du th\'eor\`eme II par les perturbations
de $\widetilde\gamma$ tangentes au feuilletage et \`a extr\'emit\'es dans $T$) est localement conjugu\'e \`a l'application $y\mapsto (y-y_0)^{-\alpha}+y_1$  qui n'est pas alg\'ebro\"\i de
en $y_0$. C'est l'{\it application de Dulac complexe de la singularit\'e $p$}.
Bien s\^ur, le th\'eor\`eme II n'est pas contredit puisque l'on a choisi un chemin
$\widetilde\gamma$ passant par une singularit\'e du feuilletage et donc relevant un chemin
$\gamma$ dans la base qui intersecte $\Sigma_E$. Cependant, toute d\'etermination locale
de $\varphi$ en un point $y_0'$ proche de $y_0$ est un germe d'application donn\'e par
le th\'eor\`eme II qui, par construction, n'admet pas de prolongement alg\'ebro\"\i de
sur $\mathbb C$.
\end{example}

\eject

\begin{figure}[htbp]
\begin{center}

\input{Alberto8.pstex_t}
 
\caption{Exemple \ref{E:Dulac}}
\label{figure:8}
\end{center}
\end{figure}

L'\'enonc\'e original du th\'eor\`eme II a d\^u suscitter maintes critiques
puisque Paul Painlev\'e a cru bon de rappeler, dans l'appendice du livre de Boutroux
(voir \cite{Bo}, p.142-144) que l'\'enonc\'e n'est que local. Pour bien montrer
qu'il avait en t\^ete qu'une interpr\'etation globale de son \'enonc\'e ne pouvait pas
\^etre vraie, il a m\^eme cru bon de donner le contre-exemple suivant. Comme nous l'avons
d\'ecortiqu\'e en d\'etails, il est plus compliqu\'e
que le pr\'ec\'edent mais tr\`es instructif pour d'autres raisons. Je remercie 
K. Okamoto de m'avoir signal\'e cet exemple en F\'evrier 2001.

\begin{example}[\cite{Bo}, p.142-144]\label{E:Painleve}\rm
Consid\'erons l'\'equation diff\'erentielle~:
$$(E)\ \ \ \ \ {dy\over dx}={y\over x(y+1)}.$$
Son int\'egration nous conduit \`a l'int\'egrale premi\`ere globale $H(x,y)={ye^y\over x}$ et
les solutions globales sont d\'efinies par
l'\'equation implicite 
$H(x,f(x,x_0,y_0))=H(x_0,y_0)$.
Apr\`es compactification dans $\overline{\mathbb C}\times\overline{\mathbb C}$, le feuilletage
$\mathcal F$ poss\`ede exactement $4$ droites invariantes, \`a savoir $\{x=0\}$, $\{x=\infty\}$,
$\{y=0\}$ et $\{y=\infty\}$. Toute feuille autre que ces droites est transcendante
et admet une param\'etrisation uniforme globale par $y\in\mathbb C$,
\`a savoir $y\mapsto c\cdot ye^y$. N\'eanmoins, les solutions correspondantes
$y=f(x,x_0,y_0)$ sont des fonctions multiformes de $x$ comme nous allons le voir.
Par ailleurs, $\mathcal F$ poss\`ede exactement $4$ points singuliers
situ\'es aux intersections de ces droites invariantes. 

\begin{figure}[htbp]
\begin{center}

\input{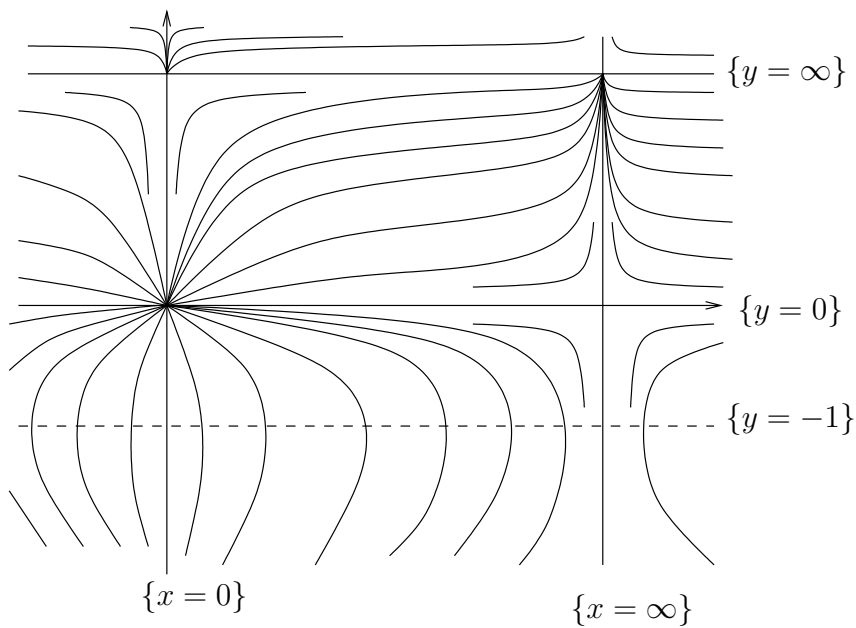}
 
\caption{Feuilletage associ\'e \`a l'exemple \ref{E:Painleve}}
\label{figure:9}
\end{center}
\end{figure}

Pr\'ecis\'ement, $(0,0)$
est un n\oe ud $x(1+y)dy-ydx=0$ et est localement conjugu\'e \`a $xdy=ydx$
d'apr\`es le th\'eor\`eme de Poincar\'e.
Ceci signifie qu'il existe un changement de coordonn\'ees analytiques locales
envoyant les feuilles de $\mathcal F$ sur celles du feuilletage radial $\{{y\over x}=constant\}$.
De m\^eme, $(\infty,0)$ est une singularit\'e selle $X(1+y)dy+ydX=0$ (ici $X={1\over x}$)
et est localement conjugu\'ee \`a sa partie lin\'eaire $Xdy+ydX=0$.
Finalement, les singularit\'es en $(0,\infty)$ et en $(\infty,\infty)$
sont des selle-n\oe uds respectivement donn\'es par $x(1+Y)dY)+Y^2dx$ (o\`u $Y={1\over y}$)
et $X(1+Y)dY)-Y^2dX$. Donc $\Sigma_E=\{0,\infty\}$.
Enfin, le lieu des tangences verticales de $\mathcal F$, en dehors des $2$ droites verticales
invariantes, est donn\'e par $\Delta=\{y=-1\}$.
La feuille $L_{x_0,y_0}$ correspondant \`a n'importe quelle solution (non constante)
$y=f(x,x_0,y_0)$
intersecte $\{y=0\}$ une fois en $(0,0)$, $\Delta$ une fois en $(x_1,-1)$ avec
$x_1=-{x_0\over y_0e^{y_0+1}}$ et accumule la droite horizontale $\{y=\infty\}$.
Ceci se d\'eduit ais\'ement de la parametrisation par la variable $y$.
Donc $f(x,x_0,y_0)$ poss\`ede exactement $3$ singularit\'es, \`a savoir en $0$, $x_1$
et $\infty$, et poss\`ede une infinit\'e de d\'eterminations en tout autre point.
Autour de $x=0$, une des branches poss\`ede une extension analytique et co\"\i ncide
avec une des feuilles locales de la singularit\'e radiale de $\mathcal F$.
Les autres branches tendent vers $\infty$
et se permutent transitivement lorsque l'on tourne autour de $x=0$. Autour de $x=\infty$,
toutes les branches tendent vers $\infty$ et se permutent transitivement.
Autour de $x_1$, seules $2$ branches
se permutent autour de la singularit\'e alg\'ebro\"\i de alors que toute autre branche
poss\`ede une extension analytique en $x_1$.
Ceci nous donne une description qualitative compl\`ete des solutions.

\begin{figure}[htbp]
\begin{center}

\input{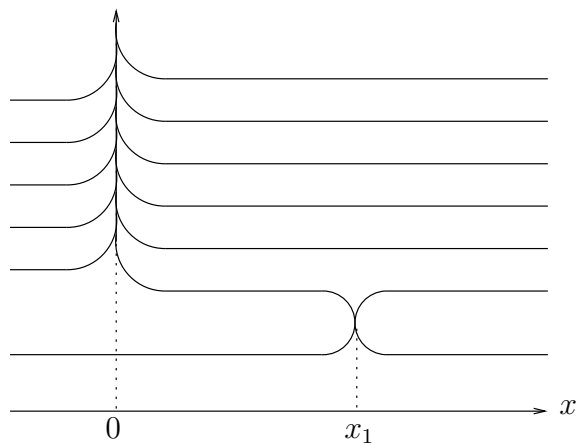}
 
\caption{Graphe d'une solution de l'exemple \ref{E:Painleve}}
\label{figure:10}
\end{center}
\end{figure}

\end{example}

\eject

\begin{remark}\rm
Painlev\'e a soulev\'e le paradoxe suivant.
Choisissons une solution analytique locale $f(x,x_0,y_0)\not\equiv0$
en un point de $\Omega$,
disons $x_0=1$ pour simplifier. En particulier, l'unique singularit\'e alg\'ebro\"\i de
$x_1$ est suppos\'ee distincte de $x_0$, c'est \`a dire $y_0\not=-1$.
Alors, toute d\'etermination de $f(x)$ en $x_0$ peut \^etre atteinte par prolongement analytique
le long d'un lacet $\gamma$ d'indice $1$ autour de $x_1$ mais d'indice $0$ autour de $0$.
En d'autres termes, une infinit\'e de branches sont permut\'ees uniquement en tournant
autour de la singularit\'e alg\'ebro\"\i de sans tourner autour des autres !
En effet, choisissons un chemin $\sigma:[0,1]\to\Omega^*$, $\Omega^*:=\Omega\setminus\{x_1\}$,
allant de $x_0$ vers un point $x_0'$ proche de $0$ le long duquel le prolongement analytique
de $f(x)$ conduit \`a l'unique branche analytique en $0$. Alors, notons $\tau:[0,1]\to\Omega^*$
un lacet proche de $0$, d'indice $1$ autour de $0$ et d'extr\'emit\'e $x_0'$.
Finalement, choisissons un lacet $\gamma_0:[0,1]\to\Omega^*$ d'extr\'emit\'e $x_0'$,
d'indice $1$ autour de $x_1$ et d'indice
$0$ autour de $0$, le long duquel le prolongement analytique de la branche analytique en $0$
est remplac\'ee par une autre d\'etermination gr\^ace \`a la singularit\'e alg\'ebro\"\i de.
Alors, toute d\'etermination de $f(x)$ en $x_0$ peut \^etre atteinte par prolongement analytique
le long de l'un des chemins $\sigma^{-1}\tau^{-n}\gamma_0\tau^n\sigma$, $n\in\mathbb Z$.
En effet, on a $f_{\tau^n\sigma}=f_{\sigma}$ puisqu'il s'agit ,de la branche analytique en $0$.
Par contre, $f_{\gamma_0\tau^n\sigma}$ est l'une des autres branches pr\`es de $0$
qui sont transitivement permut\'ees juste en tournant autour de $0$, c'est \`a dire
par prolongement analytique le long de chemins $\tau^{-n}$.
\end{remark}

\begin{example}[suite de l'exemple \ref{E:Painleve}]\rm
Fixons maintenant une verticale $T:\{x=x_0\}$, disons $x_0=1$ pour simplifier,
et consid\'erons l'application multivalu\'ee
$\varphi:T\to T;y\mapsto f(x_0,x_0,y)$ construite en recollant
tous les germes donn\'es par le th\'eor\`eme II.
Par exemple, consid\'erons $\gamma:[0,1]\to\Omega=\mathbb C^*$ un lacet d'indice $1$
autour de $0$ d'extr\'emit\'e en $x_0$. En appliquant le th\'eor\`eme II \`a $\gamma$
pour $y_0=\infty$ on d\'eduit que le germe d'application d'holonomie $\varphi_1$
de la courbe invariante $\{y=\infty\}$ du selle-n\oe ud (voir la remarque 11).
L'application multivalu\'ee doit au moins contenir ce germe ainsi que ses it\'er\'ees
$\varphi_n:=\varphi_1^{\circ n}$, $n\in\mathbb Z$, correspondant \`a $\gamma^n$.
Pr\'ecis\'ement, $\varphi_1$ est un germe de diff\'eomorphisme fixant $\infty$ de la forme
$\varphi_1(y)=y+2i\pi+\sum_{n>0}{a_n\over y^n}$.
D'autre part, en appliquant le th\'eor\`eme II \`a $\gamma$ pour $y_0=0$
on obtient le germe {\it identit\'e} puisque l'holonomie d'une singularit\'e radiale est triviale.

\eject

Bien s\^ur, le germe identit\'e et $\varphi_1$ ne peuvent pas \^etre connect\'es
par prolongement analytique le long du plan des $y$. Donc, {\it l'application multivalu\'ee $\varphi$
doit \^etre la r\'eunion de plusieurs fonctions multiformes}. En d'autres termes,
son graphe n'est pas irreducible.

Afin d'apr\'ehender l'application $\varphi:y\mapsto f(x_0,x_0,y)$ toute enti\`ere,
remarquons que deux points quelconques $y_0$ et $y_1$ sur $T$ sont dans la m\^eme feuille
si et seulement si $H(1,y_0)=H(1,y_1)$. Dans le cas o\`u $y_0,y_1\not=0$, ceci signifie encore
que $y_0+\log(y_0)=y_1+\log(y_1)$ pour des d\'eterminations convenables du logarithme.
En d'autres termes, $y_1=\varphi(y_0)$ pour une d\'etermination de $\varphi$
si et seulement si $\psi(y_0)$ et $\psi(y_1)$ diff\`erent d'un \'el\'ement de $2i\pi\mathbb Z$
(ind\'ependamment de la d\'etermination de $\psi(y)=y+\log(y)$ choisie).
Donc, $\varphi$ n'est rien d'autre que le groupe cyclique de translations $2i\pi\mathbb Z$
tir\'e en arri\`ere par la fonction multiforme $\psi$.

Introduisons le champ de vecteurs $Z=2i\pi\partial_z$ dont l'exponentielle
$\exp(Z):z\mapsto z+2i\pi$ au temps $1$ engendre le groupe cyclique pr\'ec\'edent.
Son rel\`evement par $\psi$ est $\psi^*Z=2i\pi{y\over 1+y}\partial_y$.
Sa partie r\'eelle est donn\'ee par :
$$V={2\pi\over (1+u)^2+v^2}\left((u+u^2+v^2)\partial_v-v\partial_u\right)$$
o\`u $y=u+iv$. Ainsi, $\varphi$ consiste en la collection de $\exp(nV)$, $n\in\mathbb Z$.
Ce champ de vecteurs est m\'eromorphe en $y_0=-1$ et s'annule en $0$.
Au point $y=\infty$, $\psi^*Z$ est holomorphe et $V$ est analytique,
tangent au champ de translations $2\pi\partial_v$. Donc, pr\`es de l'infini, les trajectoires
de $V$ ressemblent \`a des droites verticales.
Le portrait de phases de $V$, donn\'e par le champ de vecteurs polynomial
$( {(1+u)^2+v^2\over 2\pi})V=(u+u^2+v^2)\partial_v-v\partial_u$,
poss\`ede deux singularit\'es, \`a savoir une selle en $y_0$ et un centre en $0$.
Les $4$ s\'eparatrices de $y_0$ consistent en une trajectoire analytique $\sigma^+$
venant de l'infini depuis le bas et allant vers $y_0$,
une trajectoire $\sigma^-$ partant de $y_0$ pour aller vers l'infini par le haut
et une trajectoire cyclique $\sigma$ se refermant sur $y_0$
apr\`es avoir tourn\'e une fois autour du centre. Toute autre trajectoire vient de
l'infini par le bas pour aller vers l'infini par le haut. Le centre $0$ est isochrone pour $V$ :
$\exp(V)$ est le germe identit\'e pr\`es de $0$. Finalement, introduisons la suite de points
$y_n$ sur $\sigma^+\cup\sigma^-$, $n\in\mathbb Z^*$, telle que $\exp(nZ)$ envoit $y_n$
sur le p\^ole $y_0$. Par exemple, $y_1$ est le point limite sur $\sigma^+$ depuis lequel
on peut int\'egrer $V$ durant le temps $1$.

\eject

\begin{figure}[htbp]
\begin{center}

\input{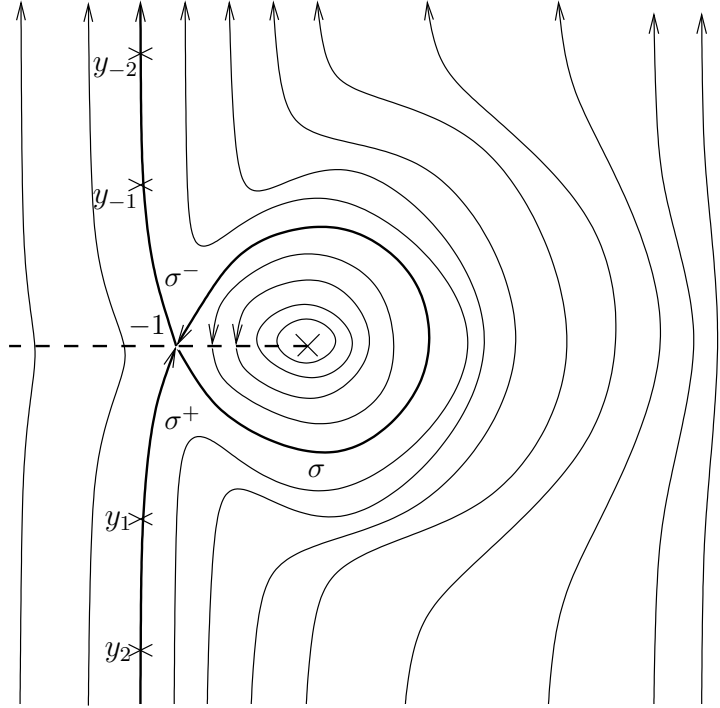}
 
\caption{Portrait de phase du champ de vecteurs (r\'eel) $V$}
\label{figure:11}
\end{center}
\end{figure}

Pour la d\'etermination principale du logarithme, la fl\`eche $\psi:y\mapsto y+\log(y)$
est un hom\'eomorphisme conforme de $\mathbb C\setminus\mathbb R^-$ sur
$\mathbb C\setminus (\mathbb R^--1\pm i\pi)$.
Pr\'ecis\'ement, $\psi$ envoit diff\'eomorphiquement la fronti\`ere :
$$\left\{\begin{matrix}
]-\infty,-1]&\text{par le bas sur}&]-\infty,-1]-i\pi&\text{par le bas,}\\
[-1,0[&\text{par le bas sur}&[-1,-\infty[-i\pi&\text{par le haut,}\\
]0,-1]&\text{par le haut sur}&]-\infty,-1]+i\pi&\text{par le bas,}\\
[-1,-\infty[&\text{par le haut sur}&[-1,-\infty[+i\pi&\text{par le haut.}\\
\end{matrix}\right.$$
Ceci se voit bien sur l'application $\psi$ \'ecrite en coordonn\'ees polaires :
$$\psi(re^{i\theta})=[r\cos(\theta)+\log(r)]+i[r\sin(\theta)+\theta].$$
Les points $y_n$ pr\'ec\'edemment d\'efinis sont les images r\'eciproques par $\psi$
des points
$z_n:=-1-(1+2n)i\pi$ pour $n>0$ et $z_n:=-1+(1-2n)i\pi$ pour $n<0$.
Finalement, la monodromie de $\psi$ est $\psi(e^{2i\pi}y)=\psi(y)+2i\pi$.

\begin{figure}[htbp]
\begin{center}

\input{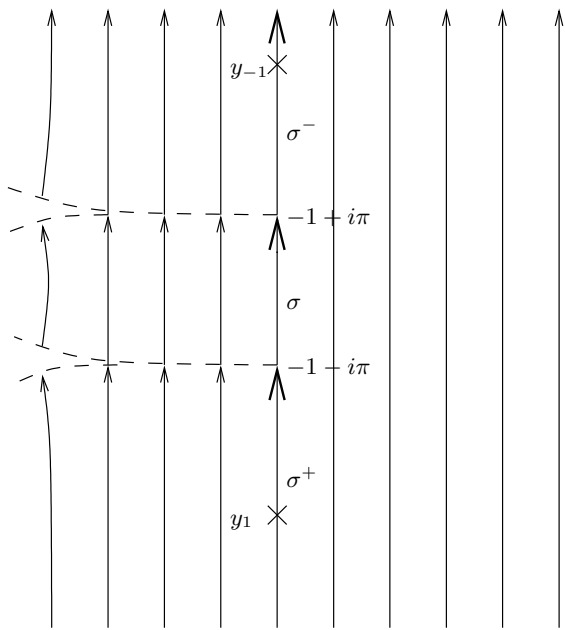}
 
\caption{Redressement de $V$ par $\psi$}
\label{figure:12}
\end{center}
\end{figure}

En raisonnant ou bien sur la dynamique de $V$ ou bien avec l'application $\psi$
telle qu'elle est d\'ecrite juste avant, on obtient la description suivante.

L'application multivalu\'ee $\varphi$ consiste en exactement deux $2$ fonctions,
\`a savoir l'application identit\'e et une fonction multiforme \`a singularit\'es
alg\'ebro\"\i des en $y_n$, $n\in\mathbb Z^*$, et $2$ singularit\'es
non alg\'ebro\"\i des en $0$ et $\infty$. En particulier, {\it toutes les it\'er\'ees
de l'holonomie $\varphi_n$, $n\in\mathbb Z^*$, sont les d\'eterminations d'une m\^eme fonction}.
Pr\'ecis\'ement, pour $n>0$, $\varphi_n$ et $\varphi_{n+1}$ sont permut\'ees autour de $y_n$
et $\varphi_{-n}$ et $\varphi_{-n-1}$ sont permut\'ees autour de $y_{-n}$. Les points
$y_n$ sont effectivement singuliers uniquement lorsqu'ils sont vus depuis ces branches.
Le point $\infty$ est r\'egulier pour toutes les d\'eterminations mais devient singulier
dans la situation suivante.
Partons d'une branche $\varphi_1$ et proc\'edons au prolongement analytique le long
d'une trajectoire $\sigma^+$ jusqu'\`a $y_1$. Alors, changeons de d\'etermination autour de
$y_1$ et revenons au point $\infty$ le long de $\sigma^+$ en tournant successivement
une fois autour de chacun des points $y_2,\ y_3,\ y_4,\ldots$ dans cet ordre.
En particulier, $\varphi$ n'a pas de limite le long de ce chemin
puisque $\varphi(y)$ parcours $\sigma$ entre deux singularit\'es successives
interm\'ediaires $y_n$, $y_{n+1}$. Par ailleurs, toutes les branches ont pour limite $\infty$
en $0$ et se permutent transitivement autour, comme pour la fonction logarithme.

Les singularit\'es alg\'ebro\"\i des $y_n$, $n\in\mathbb Z^*$, sont les points d'intersection
de la verticale $T$ avec la feuille passant par $(1,-1)$, c'est \`a dire qui est tangente \`a $T$.
Elles apparaissent lorsque, partant d'une solution analytique en $(1,y_n)$, on proc\`ede au
prolongement analytique le long d'un chemin $\gamma$ tournant $n$ fois autour de $0$
(dans le plan des $x$)
et revenant au point $1$ avec la solution alg\'ebro\"\i de en $(1,-1)$, i.e.
avec tangence sur $T$.
\end{example}

\section{Conjectures}

Le th\'eor\`eme I ainsi que les exemples pr\'ec\'edents
d\'ecrivent compl\`etement le type de singularit\'es des solutions
et la fa\c con dont elles appara\^\i ssent.
Le th\'eor\`eme II, quant \`a lui, ne nous fournit pas d'indication
sur le comportement des applications d'holonomie vis \`a vis du
prolongement analytique. L'exemple \ref{E:Dulac} nous montre d'ores et d\'ej\`a
que ces applications ne sont pas globalement alg\'ebro\"\i des~:
des singularit\'es plus compliqu\'ees peuvent appara\^\i tre.
Nous proposons deux conjectures.

\begin{conj}
Soient $P,Q\in\mathbb C[x,y]$ et consid\'erons l'\'equation diff\'erentielle :
$$(E)\ \ \ \ \ {dy\over dx}={P(x,y)\over Q(x,y)}.$$
Supposons que ni $P$, ni $Q$, ne soient identiquement nuls et consid\'erons
deux droites verticales $T_0$ et $T_1$ (non invariantes par le feuilletage)
ainsi qu'un germe d'applica\-tion alg\'ebroide
$\varphi:(T_0,p_0)\to (T_1,p_1)$ donn\'e par le Th\'eor\`eme II, $p_i\in T_i$, $i=0,1$.
Alors l'ensemble $\Sigma_\varphi$ des singularit\'es de $\varphi$
pour le prolongement analytique dans $T_0\simeq\overline{\mathbb C}$ est au plus d\'enombrable.
\end{conj}

En particulier, $\varphi$ admet un prolongement sans coupure :
le ph\'enom\`eme de l'exemple \ref{E:Picard} ne se produit pas.
Cette conjecture fut motiv\'ee au d\'epart par le probl\`eme suivant.
Si l'\'enonc\'e est vrai, alors on peut construire une int\'egrale
premi\`ere multiforme pour le feuilletage dont l'ensemble singulier est une
r\'eunion d\'enombrable de feuilles. Une telle fonction admet automatiquement
un groupe de monodromie (voir \cite{Kh}) qui sera en g\'en\'eral 
tr\`es consistant.
\`A ma connaissance, on n'a jamais associ\'e que des pseudo-groupes
\`a un feuilletage alg\'ebrique g\'en\'eral, ce qui rend difficile
la construction d'une th\'eorie de Galois pour ce type d'objets.

\eject

Dans le cas o\`u les singularit\'es de l'\'equation (E) sont toutes hyperboliques,
il est vraisemblable que  les singularit\'es non alg\'ebro\"\i des de $\varphi$
naissent des applications de Dulac locales des singularit\'es. On peut alors
s'attendre \`a ce qu'elles se situent sur les points d'intersection de la droite
$T_0$ avec les feuilles (en nombre fini) ``passant par ces singularit\'es''
sous forme de courbes invariantes locales.
Notamment, l'application $\varphi$, apr\`es prolongement analytique maximal,
mettrait en relation des points de $T_0$ et $T_1$ qui ou bien seraient
dans la m\^eme feuille, ou bien seraient dans des feuilles distinctes
se connectant comme courbes invariantes locales aux singularit\'es
(comme un polycycle).
C'est en tout cas ce qui se passe dans tous les exemples \'etudi\'es.
Le lecteur pourra trouver, dans la th\`ese de David Marin (voir \cite{Ma}),
d'autres exemples instructifs illustrant cette conjecture.

Cependant, dans tous les exemples \'etudi\'es, on utilise une int\'egrale premi\`ere
multiforme explicite du feuilletage (que l'on inverse) pour pouvoir obtenir une description
globale de $\varphi$. Nous ne sommes pas capable actuellement de donner un seul exemple
d'\'equation non int\'egrable (explicitement) satisfaisant notre conjecture.
Ces applications $\varphi$ peuvent \^etre extr\`emement compliqu\'ees.
Par exemple, un travail r\'ecent (voir \cite{BeLiLo}) montre que ces applications
ne satisfont g\'en\'eralement aucune \'equation diff\'erentielle analytique et sont, en ce sens,
plus transcendantes que toutes les fonctions solutions de toutes les \'equations
diff\'erentielles. Dans cette direction, nous proposons cette autre :

\begin{conj}
Soient $P,Q\in\mathbb C[x,y]$ et consid\'erons l'\'equation diff\'erentielle :
$$(E)\ \ \ \ \ {dy\over dx}={P(x,y)\over Q(x,y)}.$$
Supposons que ni $P$, ni $Q$, ne soient identiquement nuls et que $Q$
ne contient aucun facteur vertical :
aucune droite verticale n'est totalement tangente au feuilletage.
Consid\'erons alors une droite verticale $T_0$ 
et deux germes d'applications analytiques
$\varphi:(T_0,p_0)\to (T_0,p_1)$ et 
$\varphi':(T_0,p_0)\to (T_0,p_1')$ donn\'es par le Th\'eor\`eme II. 
Alors $\varphi'$ se d\'eduit par prolongement analytique
de $\varphi$ le long de $T_0$... sauf si l'un des deux germes est l'identit\'e.
\end{conj}

\eject

Une cons\'equence imm\'ediate de la conjecture 2 et de \cite{BeLiLo} serait :

{\it si $P$ et $Q$ sont suffisamment g\'en\'eriques de degr\'e $d\ge2$, alors le prolongement analytique maximal de $\varphi$ 
accumule tout le pseudo-groupe des transformations conformes
sur $T_0$ : pour tout germe d'application analytique inversible
$\psi:(T_0,p_0)\to T_1$, il existe une suite de d\'eterminations
$\varphi_n:(T_0,p_0)\to T_1$ de $\varphi$ en $p_0$
convergeant uniform\'ement vers $\psi$ sur un voisinage de $p_0$
dans $T_0$.}

Il serait joli de construire un tel exemple.

\backmatter

\end{document}